\documentclass[11pt,leqno]{amsart}
\usepackage{amsmath,amssymb,amsfonts,amstext,amsthm}
\usepackage{mathtools,latexsym,amscd}
\usepackage[T1]{fontenc}
\usepackage[utf8]{inputenc}
\usepackage{lmodern}
\usepackage{microtype}
\usepackage{color}
\usepackage{parskip}
\usepackage{enumitem}
\usepackage{tikz}
\usetikzlibrary{arrows.meta}
\usepackage[pagebackref,colorlinks]{hyperref}
\allowdisplaybreaks[1]
\setlist[enumerate]{leftmargin=2em,itemsep=3pt,topsep=5pt}
\theoremstyle{plain}
\newtheorem{theorem}{Theorem}
\newtheorem{lemma}[theorem]{Lemma}

\newtheorem{proposition}[theorem]{Proposition}
\newtheorem{corollary}[theorem]{Corollary}
\theoremstyle{definition}
\newtheorem{remark}[theorem]{Remark}
\newtheorem*{remark*}{Remark}
\newtheorem*{question*}{Question}

\newcommand{\R}{\mathbb R}
\newcommand{\Z}{\mathbb Z}
\newcommand{\F}{F_2}
\newcommand{\Q}{\operatorname Q}
\newcommand{\Hb}{H_b}
\newcommand{\rHb}{\overline H_b}
\newcommand{\rEH}{\overline{EH}_b}
\newcommand{\Sym}{\operatorname{Sym}}
\newcommand{\Homeo}{\operatorname{Homeo}}
\newcommand{\Diff}{\operatorname{Diff}}
\newcommand{\Ham}{\operatorname{Ham}}
\newcommand{\Mod}{\operatorname{Mod}}
\newcommand{\id}{\operatorname{id}}
\newcommand{\im}{\operatorname{im}}

\newcommand{\eps}{\varepsilon}

\hypersetup{
 pdftitle={Quasimorphisms and bounded cohomology of transformation groups of nonorientable surfaces},
 pdfauthor={Michael Brandenbursky},
 pdfsubject={Gambaudo--Ghys constructions, undistortion and bounded cohomology},
 pdfkeywords={quasimorphisms, bounded cohomology, nonorientable surfaces, surface braid groups}
}
\begin{document}
\raggedbottom
\title[Quasimorphisms and bounded cohomology]
{Quasimorphisms and bounded cohomology of transformation groups of nonorientable surfaces}
\author{Michael Brandenbursky}
\thanks{M.B. was partially supported by the Israel Science Foundation grant 823/23.}
\address{Department of Mathematics, Ben Gurion University, Israel}
\email{brandens@bgu.ac.il}
\date{}
\begin{abstract}
Let $N_g$ be a closed connected nonorientable surface of genus
$g\ge1$, let $\omega$ be a smooth positive density, and let $\mu$
be a good measure. In this paper we generalize the Gambaudo--Ghys
construction using homogeneous quasimorphisms on the central
quotients of the braid groups of $N_g$. We prove injectivity of this map for
density-preserving diffeomorphisms. For $g\ge2$ and quasimorphisms
vanishing on the disc braid group, we obtain an injective map to
$C^0$-continuous homogeneous quasimorphisms on
$\Homeo_0(N_g,\mu)$. We also obtain undistortion criteria in
$\Homeo_0(N_g)$ and $\Diff_0(N_g)$.
For $g=1,2$ and every $d\ge2$, we prove 
$$
\dim\rHb^d(\Homeo_0(N_g,\mu))=\infty,\qquad
\dim\rEH^d(\Diff_0(N_g,\omega))=\infty.
$$
Moreover, $\rEH^d(\Homeo_0(N_2,\mu))$ is infinite dimensional in every positive even
degree $d$. Here $\rHb^d$ and $\rEH^d$ denote reduced bounded cohomology
and reduced exact bounded cohomology, respectively.
\end{abstract}
\maketitle
\markboth{Quasimorphisms and bounded cohomology}
{Quasimorphisms and bounded cohomology}

\section{Introduction}\label{sec:introduction}
Quasimorphisms are useful tools for studying the algebraic
structure of nonabelian groups, especially when the space of
real-valued homomorphisms is finite dimensional; see, e.g.,
\cite{Calegari2009,EntovPolterovich2003,BestvinaFujiwara2002}.
The constructions of Gambaudo--Ghys and Polterovich produce
quasimorphisms on groups of measure-preserving transformations 
of compact connected orientable manifolds from
quasimorphisms on braid groups and fundamental groups
\cite{GambaudoGhys2004,Polterovich2006}. These constructions have been
used to study invariant metrics and dynamics of transformation groups;
see, for example, \cite{BrandenburskyKedra2013,Brandenbursky2015,Entropy}.

In \cite{BKMS2026}, Brandenbursky, K\k edra, Marcinkowski and Shelukhin
identified the Gambaudo--Ghys quasimorphisms which extend continuously
to measure-preserving homeomorphisms of closed orientable surfaces of
positive genus, and obtained applications to distortion in the full
homeomorphism group.
In this paper we study the corresponding construction on closed
nonorientable surfaces and its applications to bounded cohomology.
The new construction is used to obtain infinite-dimensional reduced
exact bounded cohomology in every even degree for the
measure-preserving homeomorphism group of the Klein bottle.

Let $N_g$ be the closed nonorientable surface of genus $g$. Thus
$N_1=\mathbb{RP}^2$ and $N_2$ is the Klein bottle. A smooth positive
density $\omega$ on $N_g$ defines a measure $\mu$. On an oriented open
set, this density is represented by a positive area form. Let
$\Diff_0(N_g,\omega)$ and $\Homeo_0(N_g,\mu)$ denote the identity
components of the density-preserving diffeomorphism group and the
measure-preserving homeomorphism group, respectively. A finite Borel
measure is called \emph{good} if it is nonatomic and positive on every
nonempty open set. We normalize the measures defining the transformation groups
to have total mass one. All group-cohomology coefficients are trivial and real, and
transformation groups are regarded as discrete when taking cohomology.

\emph{The Gambaudo--Ghys construction.}
Let $B_n(N_g)$ be the braid group and put
$$
\widehat B_n(N_g)=B_n(N_g)/Z(B_n(N_g)),
$$
where $Z(B_n(N_g))$ denotes the center of $B_n(N_g)$. 
For a group $K$, denote by $\Q(K)$ its vector space of homogeneous
quasimorphisms. Let $\gamma_n$ be the trajectory cocycle and $X_n(N_g)$ the
ordered configuration space. For $q\in\Q(\widehat B_n(N_g))$, the
Gambaudo--Ghys construction is the homogenization of the map
$$
f\longmapsto\int_{X_n(N_g)}q(\gamma_n(f,x))\,d\mu^n(x).
$$
The central quotient removes the dependence on the isotopy used to
define the trajectory braid. The definitions are given
in Section~\ref{sec:preliminaries}.

\begin{theorem}\label{thm:smooth-main}
Let $g\ge1$, $n\ge2$, and let $\omega$ be a smooth positive density
on $N_g$. The Gambaudo--Ghys construction defines an injective linear map
$$
\mathcal G_n^\omega:\Q(\widehat B_n(N_g))\longrightarrow\Q(\Diff_0(N_g,\omega)).
$$
It is also injective after quotienting both spaces by their
homomorphisms. Its image modulo homomorphisms is infinite dimensional
if $g\ge2$, or if $g=1$ and $n\ge3$.
\end{theorem}

For $g=1$ and $n=2$, the braid group $B_2(N_1)$ is finite,
and hence $\Q(\widehat B_2(N_1))=0$. For the construction on
homeomorphisms, we assume $g>1$ and impose the following
vanishing condition. Fix an embedded disc $D\subset N_g$
containing the base configuration, and let
$\iota_D:B_n(D)\longrightarrow\widehat B_n(N_g)$
be the homomorphism induced by the inclusion $D\hookrightarrow N_g$
followed by the quotient map
$B_n(N_g)\to\widehat B_n(N_g)$.
The map $\iota_D$ is injective by
\cite[Propositions~1.6 and~2.2]{ParisRolfsen1999}
and the description of the center of the Klein-bottle braid
group in \cite[Proposition~5.2]{GuaschiPereiro2018}.
We define
$$
\Q_{\mathrm{disc}}(\widehat B_n(N_g))=
\bigl\{q\in\Q(\widehat B_n(N_g)):q\circ\iota_D=0\bigr\}.
$$
Write $\Q_{C^0}(K)$ for the space of continuous homogeneous
quasimorphisms on a group $K$ equipped with the $C^0$ topology.

\begin{theorem}\label{thm:homeo-main}
Let $g\ge2$, $n\ge2$, and let $\mu$ be a good measure on $N_g$.
The Gambaudo--Ghys construction defines an injective linear map
$$
\mathcal G_n^\mu:\Q_{\mathrm{disc}}(\widehat B_n(N_g))\longrightarrow\Q_{C^0}(\Homeo_0(N_g,\mu)).
$$
It is also injective modulo homomorphisms, and its image modulo
homomorphisms is infinite dimensional for every fixed $n\ge2$.
For a smooth density, it extends the map of
Theorem~\ref{thm:smooth-main} on this subspace.
\end{theorem}

For $g\ge3$, an infinite-dimensional space of homogeneous
quasimorphisms on $\Homeo_0(N_g,\mu)$ is already obtained
in the one-point case, starting from homogeneous
quasimorphisms on $B_1(N_g)=\pi_1(N_g)$.
This infinite-dimensionality was proved by
Kumar~\cite{Kumar2021}; see also
Brandenbursky--Menashe~\cite{BrandenburskyMenashe2026}.
Theorem~\ref{thm:homeo-main} establishes injectivity and
infinite-dimensional image for every fixed number of
strands $n\ge2$. 

For the Klein bottle, every homogeneous quasimorphism on
$\pi_1(N_2)$ is a homomorphism, since this group is amenable.
We obtain the required infinite-dimensional families from
homogeneous quasimorphisms on $\widehat B_n(N_2)$,
$n\ge2$, which vanish on the disc braid subgroup; see
Proposition~\ref{prop:disc-vanishing-abundance}.
The injectivity argument applies to all quasimorphisms in the indicated
spaces. It uses invariant smooth positive densities concentrating most
of their mass in small discs, and a polynomial whose mixed coefficient
detects a prescribed braid.

Kim and Maruyama constructed an injective two-strand Gambaudo--Ghys
map for density-preserving diffeomorphisms of the M\"obius band
\cite{KimMaruyamaIMRN2026}. Here we treat every number of strands on
closed nonorientable surfaces and establish injectivity and continuous
extension to the measure-preserving homeomorphism groups. For the smooth
construction we first establish a finite-braid lemma by an autonomous
flux correction and fragmentation into Hamiltonian disc maps, following
the argument of \cite[Section~2.B.2]{Entropy}.

\emph{Undistortion.}
Write $\Homeo_0(N_g)$ and $\Diff_0(N_g)$ for the identity components
without a measure-preservation condition. An element is undistorted
if its powers have positive stable word length in every finitely
generated subgroup containing it.
Fix a smooth positive probability density $\omega$, with associated
measure $\mu$. In the next theorem this is only a reference measure
for integration; the transformations need not preserve it.

\begin{theorem}\label{thm:ambient-cochains}
Let $n\ge2$. For $g\ge2$ and
$q\in\Q_{\mathrm{disc}}(\widehat B_n(N_g))$, the formula
$$
 \Psi_{n,q}^{\omega}(f)
 =\int_{X_n(N_g)}q(\gamma_n(f,x))\,d\mu^n(x)
$$
defines a cochain on $\Homeo_0(N_g)$ with semi-bounded differential:
$$
 \sup_{h\in\Homeo_0(N_g)}
 \left|\Psi_{n,q}^{\omega}(f)+\Psi_{n,q}^{\omega}(h)
                   -\Psi_{n,q}^{\omega}(fh)\right|<\infty
 \qquad(f\in\Homeo_0(N_g)).
$$
For $g\ge1$ and arbitrary $q\in\Q(\widehat B_n(N_g))$, the same
assertion holds on $\Diff_0(N_g)$.
In either case, if
$$
 \limsup_{k\to\infty}\frac{|\Psi_{n,q}^{\omega}(f^k)|}{k}>0,
$$
then $f$ is undistorted in the corresponding ambient group.
\end{theorem}

On the measure-preserving subgroup, the homogenization is the
quasimorphism constructed above. On the ambient group, we assert
semi-boundedness, not the quasimorphism property. Pointwise evaluation
of the trajectory cocycle gives a second criterion, which does not
require a nonzero averaged value.
Suppose that $x\in X_n(N_g)$ is fixed by $f^r$ for some $r\ge1$.
The trajectories of this configuration define a return braid
$\beta_{f,x,r}\in\widehat B_n(N_g)$, well defined up to conjugation.

\begin{theorem}\label{thm:undistortion-main}
Let $n\ge2$.
\begin{enumerate}[label=\textup{(\roman*)}]
\item
Let $g\ge2$ and $f\in\Homeo_0(N_g)$.
If an $n$-point periodic configuration of $f$ has return braid $\beta$
with $q(\beta)\ne0$ for some
$q\in\Q_{\mathrm{disc}}(\widehat B_n(N_g))$, then $f$ is
undistorted in $\Homeo_0(N_g)$.
\item
Let $g\ge1$ and $f\in\Diff_0(N_g)$.
If an $n$-point periodic configuration of $f$ has return braid $\beta$
with $q(\beta)\ne0$ for some $q\in\Q(\widehat B_n(N_g))$, then
$f$ is undistorted in $\Diff_0(N_g)$.
\end{enumerate}
\end{theorem}

These are the nonorientable versions of the semi-bounded and
pointwise constructions in
\cite[Theorem~1.1, Corollaries~1.7--1.8 and Proposition~4.1]{BKMS2026}.
For diffeomorphisms we also prove that an infinite-order return braid
suffices; see Theorem~\ref{thm:smooth-undistortion}.
The multiple-point assertions for arbitrary homeomorphisms are not
consequences of \cite[Theorem~1.3]{FranksHandel2006}, which concerns
$C^1$ diffeomorphisms. That theorem already treats the full
$C^1$ identity component, so removing measure preservation alone is
not a novelty claim in the smooth category. We discuss the overlap
with the displacement and smooth dynamical obstructions of
Franks--Handel at the end of Section~\ref{sec:undistortion}.

\emph{Bounded cohomology.}
Denote bounded cohomology by $\Hb^d$, its reduced version by $\rHb^d$,
and reduced exact bounded cohomology by $\rEH^d$.
Brandenbursky--Marcinkowski \cite{BrandenburskyMarcinkowski2022},
Kimura \cite{Kimura2025}, Nitsche \cite{Nitsche}, and
Brandenbursky--Menashe \cite{BrandenburskyMenashe2026} used transfer
maps to construct higher-degree bounded classes on transformation
groups. In \cite{Brandenbursky2026} the author 
used product and mixed transfers in order to obtain
infinite-dimensional families in every degree greater than one for
orientable surfaces. The nonorientable hyperbolic case was treated in
\cite{BrandenburskyNonorientable2026}. Here we address the two remaining
closed nonorientable surfaces. Their fundamental groups are amenable,
and hence have vanishing bounded cohomology in positive degrees
\cite[Chapter~3]{Frigerio}. Thus transfer from the fundamental group
cannot provide the nonzero classes needed in these cases.

\begin{samepage}
\begin{theorem}\label{thm:bounded-main}
Let $g\in\{1,2\}$, let $\mu$ be a good measure on $N_g$, and let
$\omega$ be a smooth positive density. For every $d\ge2$,
$$
\dim\rHb^d(\Homeo_0(N_g,\mu))=\infty,\qquad\dim\rEH^d(\Diff_0(N_g,\omega))=\infty.
$$
Moreover, for every $k\ge1$,
$$
\dim\rEH^{2k}(\Homeo_0(N_2,\mu))=\infty.
$$
\end{theorem}
\end{samepage}

The corresponding unreduced cohomology groups are therefore infinite
dimensional as well. Together with
\cite{BrandenburskyNonorientable2026}, the theorem gives
infinite-dimensional reduced bounded cohomology of
$\Homeo_0(N_g,\mu)$ in every degree $d\ge2$ for every $g\ge1$.
For density-preserving diffeomorphisms, the classes are exact in every
genus. For homeomorphisms, the exact assertions in Theorem~\ref{thm:bounded-main}
are stated separately.

For $N_1=\mathbb{RP}^2$, we use the virtually free quotient of the
three-strand braid group; see Lemma~\ref{lem:rp2-source}.
For the Klein bottle, we use exact bounded classes on
$\widehat B_2(N_2)$ whose restrictions to a two-generator free
subgroup have infinite-dimensional image in degrees two and three;
see Lemma~\ref{lem:klein-source}. This subgroup is represented by
moving one point along orientation-preserving loops while fixing
the other point.
In Section~\ref{sec:bounded-applications}, we construct
density-preserving pushing maps which detect the resulting
degree-two and degree-three classes after transfer.
The pushing maps are made density preserving by applying
the relative Moser argument on the orientation double cover.
The correcting isotopies commute with the deck involution
and therefore descend to the Klein bottle.
The degree-three classes provide the odd-degree part of
Theorem~\ref{thm:bounded-main}.

We combine these classes with fixed degree-two braid classes
using the injective external-product and mixed-transfer
constructions of \cite{Brandenbursky2026}.
The new Gambaudo--Ghys construction enters the exactness
argument for Klein-bottle homeomorphisms: the exact degree-two
classes on $\widehat B_2(N_2)$ are represented by homogeneous
quasimorphisms which vanish on the disc braid group.
Theorem~\ref{thm:homeo-main} therefore supplies homogeneous
quasimorphisms on $\Homeo_0(N_2,\mu)$.
Their bounded coboundaries make the two-point factor exact,
yielding the even-degree exact assertion.
For diffeomorphisms, the finite-braid lemma allows ordinary
cochain transfer and proves exactness in every degree.

Section~\ref{sec:preliminaries} fixes the conventions and
recalls the transfer maps.
Sections~\ref{sec:smooth-gg} and \ref{sec:homeomorphisms}
prove Theorems~\ref{thm:smooth-main} and~\ref{thm:homeo-main}.
Section~\ref{sec:undistortion} proves the semi-bounded extension and
undistortion results in Theorems~\ref{thm:ambient-cochains}
and~\ref{thm:undistortion-main}.
Section~\ref{sec:low-genus} establishes the braid-group
lemmas for $N_1$ and $N_2$, and
Section~\ref{sec:bounded-applications} proves
Theorem~\ref{thm:bounded-main}.

\textbf{Acknowledgments.}
The author acknowledges the support of Israel Science Foundation grant
823/23. He thanks the Department of Mathematics at Ben Gurion University for
excellent working conditions. The author used an LLM for computational exploration, algebraic checks,
and editorial assistance. The mathematical ideas, research direction, Gambaudo-Ghys embedding and the transfer method
approach were developed exclusively by the author.

\section[Preliminaries]{Preliminaries: densities, braid groups, trajectory cocycle, quasimorphisms and bounded cohomology}
\label{sec:preliminaries}

We write $N_g$ for the closed connected nonorientable surface of
genus $g$, and $\Sigma_g$ for the closed oriented surface of genus $g$.
We write $N_{g,1}$ for $N_g$ with the interior of a disc removed.

\emph{Densities and good measures.}
A smooth positive density on $N_g$ is a positive smooth section of
$|\Lambda^2T^*N_g|$. It determines a finite Borel measure, which we normalize
to have total mass one. We denote the density by $\omega$ and the resulting
measure by $\mu$. On an oriented open subset $U\subset N_g$, the restriction
of $\omega$ is represented by a unique positive area form $\omega_U$.
For an oriented disc $D$, let $\Ham_c(D)$ denote the group of compactly
supported Hamiltonian diffeomorphisms with respect to this area form,
and let $\operatorname{Cal}_D$ denote its Calabi homomorphism.
The orientation double cover
$$
p:\Sigma_{g-1}\longrightarrow N_g
$$
has an area form $\widetilde\omega$ such that
$$
\tau^*\widetilde\omega=-\widetilde\omega,\qquad
\mu(A)=\frac12\int_{p^{-1}(A)}\widetilde\omega,
\qquad \int_{\Sigma_{g-1}}\widetilde\omega=2,
$$
where $\tau$ is the nontrivial deck transformation and $A\subset N_g$
is Borel.
This is the convention used in \cite{BrandenburskyMenashe2026}.
The canonical orientation-preserving lift of a density-preserving
diffeomorphism preserves $\widetilde\omega$.

A finite Borel measure is called \emph{good} if it is nonatomic and positive
on every nonempty open set. Multiplying a good measure by a positive constant
does not change its group of measure-preserving homeomorphisms. We therefore
normalize good measures as well.

\begin{lemma}\label{lem:good-measures}
Let $\mu$ and $\nu$ be good probability measures on $N_g$. There is a
homeomorphism $h$ such that $h_*\mu=\nu$. Conjugation by $h$ induces a
topological group isomorphism
$$
h\Homeo_0(N_g,\mu)h^{-1}=\Homeo_0(N_g,\nu).
$$
This isomorphism induces isometric isomorphisms on reduced bounded cohomology and on
reduced exact bounded cohomology.
\end{lemma}
\begin{proof}
The first assertion is the Oxtoby--Ulam theorem
\cite[Part~II, Corollary~1]{OxtobyUlam1941}; see also
\cite[Theorem~3.1]{Fathi1980}. The conjugation isomorphism and its inverse
give inverse isometries of bounded cochain complexes. They commute with
the ordinary comparison maps, which proves the last assertion.
\end{proof}

Unless stated otherwise, $\mu$ will be induced by a fixed smooth positive
probability density $\omega$. Results for an arbitrary good measure follow
from Lemma~\ref{lem:good-measures}. We denote by $\Diff_0(N_g,\omega)$ and $\Homeo_0(N_g,\mu)$
the identity components of the groups of smooth
$\omega$-preserving diffeomorphisms and $\mu$-preserving
homeomorphisms, taken in the $C^\infty$ and compact-open
topologies, respectively. When discussing continuity,
we equip both groups with the compact-open ($C^0$) topology.
We write $\Diff_0(N_g)$ and $\Homeo_0(N_g)$ for the corresponding
identity components without the density- or measure-preserving
condition. Their definitions use the $C^\infty$ and compact-open
topologies, respectively.

\emph{Braid groups.}
For $n\ge1$, let
$$
X_n(N_g)=\{(x_1,\ldots,x_n)\in N_g^n:x_i\ne x_j\text{ for }i\ne j\}
$$
be the ordered configuration space. Choose a base configuration
$z=(z_1,\ldots,z_n)$ and put
$$
P_n(N_g)=\pi_1(X_n(N_g),z),\qquad
B_n(N_g)=\pi_1(X_n(N_g)/\mathfrak S_n,[z]).
$$
Here $\mathfrak S_n$ acts by permuting the coordinates.
The product measure $\mu^n$ gives $X_n(N_g)$ total measure one.
We use the central quotient
$$
\widehat B_n(N_g)=B_n(N_g)/Z(B_n(N_g)).
$$
For $n\ge2$, the centers of the pure and full surface braid groups agree in
the cases considered here. For $g\ge3$ the center is trivial; for the Klein
bottle it is infinite cyclic; for $N_1=\mathbb{RP}^2$ it has order two.
We write $\widehat P_n(N_g)$ for the image of $P_n(N_g)$ in
$\widehat B_n(N_g)$.

We denote by $\Mod(N_g,n)$ the mapping class group of $N_g$
with $n$ marked points, allowing permutations of the marked points.
For an oriented surface, $\Mod^+(\Sigma_g,n)$ denotes the corresponding
group of orientation-preserving mapping classes.
For $n\ge2$, point-pushing gives an injective homomorphism
$$
\widehat B_n(N_g)\longrightarrow\Mod(N_g,n)
$$
whose image is the kernel of forgetting the marked points. For $g\ge3$
and $g=1$, respectively, see
\cite[equations~(1.4) and~(3.3)]{GoncalvesGuaschiMaldonado2018}.
For $g=2$, use equation~(3.6) and Lemma~4.1 of the same paper,
together with \cite[Proposition~5.2]{GuaschiPereiro2018}.

Let $D\subset N_g$ be an embedded disc containing the base
configuration. For $g>1$, the inclusion $D\hookrightarrow N_g$
induces an injective homomorphism
$B_n(D)\hookrightarrow B_n(N_g)$; see
\cite[Proposition~2.2]{ParisRolfsen1999}.
For $g>2$, the center of $B_n(N_g)$ is trivial by
\cite[Proposition~1.6]{ParisRolfsen1999}.
For $g=2$, the disc braid subgroup has trivial intersection
with the center of $B_n(N_2)$, as follows from the description
of this center in
\cite[Proposition~5.2]{GuaschiPereiro2018}.
Consequently, for every $g>1$, inclusion followed by the
quotient map gives an injective homomorphism
$$
\iota_D:B_n(D)\hookrightarrow\widehat B_n(N_g).
$$
We identify $B_n(D)$ with its image under $\iota_D$.

\emph{The trajectory cocycle.}
We use the inverse-trajectory convention of
\cite{Brandenbursky2026,BrandenburskyNonorientable2026}.
Fix a constant-curvature metric on $N_g$ and choose a polygonal
fundamental domain in its universal cover such that the image $\Delta\subset N_g$
of its interior has full measure and contains all points
of the base configuration $z$. The standard short geodesic paths in the
polygon give, for almost every $x\in X_n(N_g)$, a path
$\alpha_{x}$ from $z$ to $x$ in $X_n(N_g)$.
The configurations for which two strands of such a closing path meet form
a measure-zero set. Choose the remaining paths measurably as specified
in Section~\ref{sec:undistortion}; there we make a controlled choice
at every configuration. This is the closing-path construction of
\cite[Section~2]{BKMS2026}, completed on its exceptional set.

For $f\in\Homeo_0(N_g)$, choose an isotopy $f_t$ from the identity
to $f$, and define
\begin{equation}\label{eq:trajectory}
\gamma_n(f,x)=\left[\alpha_{x}*\{f_t(x)\}*
\alpha_{fx}^{-1}\right]^{-1}\in\widehat P_n(N_g)\subset\widehat B_n(N_g).
\end{equation}
Here $fx=(f(x_1),\ldots,f(x_n))$, paths are concatenated from left
to right, and multiplication of homeomorphisms is composition, $fh=f\circ h$.
The evaluation of a loop of homeomorphisms gives a central braid.
Consequently, \eqref{eq:trajectory} is independent of the chosen isotopy
after passage to the central quotient. It satisfies
\begin{equation}\label{eq:cocycle}
\gamma_n(fh,x)=\gamma_n(f,hx)\gamma_n(h,x)
\end{equation}
at every configuration once the closing paths are fixed everywhere.
For measure-preserving transfers only the almost-everywhere values
are relevant. On the ambient groups, Section~\ref{sec:undistortion}
uses this cocycle to construct cochains with semi-bounded differentials;
no measure-preserving transfer is asserted there.

\emph{Quasimorphisms.}
Let $K$ be a group. A function $q:K\to\R$ is a quasimorphism if
$$
D(q):=\sup_{u,v\in K}|q(uv)-q(u)-q(v)|<\infty.
$$
The number $D(q)$ is called its defect. A quasimorphism is homogeneous if
$q(u^k)=kq(u)$ for all $u\in K$ and $k\in\Z$. Every quasimorphism has a
homogenization
$$
\overline q(u)=\lim_{k\to\infty}\frac{q(u^k)}k.
$$
Homogeneous quasimorphisms are invariant under conjugation and additive
on commuting elements \cite{Calegari2009}. We denote their vector space by $\Q(K)$.

For $g\ge2$ and $n\ge2$, put
\begin{equation}\label{eq:disc-space}
\Q_{\mathrm{disc}}(\widehat B_n(N_g)):=\{q\in\Q(\widehat B_n(N_g)):q|_{\iota_D(B_n(D))}=0\}.
\end{equation}
The vanishing condition is independent of the embedded disc containing the
base configuration. Indeed, the corresponding disc braid subgroups are
conjugate in the full surface braid group. When $f$ fixes every point of
$x$, relabeling the configuration conjugates its trajectory braid
in the full braid group. Consequently $q(\gamma_n(f,x))$ is
unchanged by this relabeling. 

For a measure-preserving map $f$, whenever the integral is defined, write
\begin{equation}\label{eq:gg}
\Psi_{n,q}^{\omega}(f)=\int_{X_n(N_g)}q(\gamma_n(f,x))\,d\mu^n(x),\qquad
\mathcal G_n^{\omega}(q)(f)=\lim_{k\to\infty}\frac{\Psi_{n,q}^{\omega}(f^k)}k.
\end{equation}
We prove the required integrability below. For measure-preserving maps,
\eqref{eq:cocycle} and change of variables give
$D(\Psi_{n,q}^{\omega})\le D(q)$. Thus the second expression in
\eqref{eq:gg} is the usual homogenization.

\emph{Bounded cohomology.}
For a discrete group $K$, let $C_b^d(K)^K$ be the space of bounded functions
$c:K^{d+1}\to\R$ satisfying
$$
c(k_0h,\ldots,k_dh)=c(k_0,\ldots,k_d).
$$
We use the homogeneous differential
$$
(\delta c)(k_0,\ldots,k_{d+1})=\sum_{i=0}^{d+1}(-1)^i c(k_0,\ldots,\widehat k_i,\ldots,k_{d+1}).
$$
Its cohomology is $\Hb^d(K)$. The supremum norm on cochains induces the
canonical seminorm on bounded cohomology. We write
$$
\rHb^d(K)=\Hb^d(K)/\{u:\|u\|=0\}
$$
for reduced bounded cohomology. Dropping boundedness gives the ordinary
cochain complex and the comparison map $c_K^d:\Hb^d(K)\to H^d(K)$.
Following \cite{FrigerioPozzettiSisto2015}, put
$$
 EH_b^d(K)=\ker c_K^d,\qquad
 \rEH^d(K)=\im\bigl(EH_b^d(K)\longrightarrow\rHb^d(K)\bigr).
$$
Thus reduced exact bounded cohomology is exact bounded cohomology modulo
its zero-seminorm subspace. We do not define it as the kernel of a comparison
map on reduced bounded cohomology.

In degree two, the coboundary of $q\in\Q(K)$ is represented in these
conventions by
$$
(\delta q)(k_0,k_1,k_2)=q(k_0k_1^{-1})+q(k_1k_2^{-1})-q(k_0k_2^{-1}).
$$
Every exact bounded class in degree two arises in this way, and
$[\delta q]=0$ if and only if $q$ is a homomorphism. There are no nonzero
degree-two bounded cohomology classes of zero seminorm; see
\cite{MatsumotoMorita1985,Frigerio}. In particular, a nontrivial exact
degree-two class remains nontrivial after reduction.

Pullback by a homomorphism is seminorm nonincreasing and commutes with
comparison. Quotienting a group by an amenable normal subgroup induces an
isometric isomorphism in bounded cohomology. We use these facts both before
and after reduction; see \cite{Frigerio}.

\emph{Configuration-space transfers.}
Let $\mathcal H=\Homeo_0(N_g,\mu)$. For a bounded homogeneous cochain on
$\widehat B_n(N_g)$, define
\begin{equation}\label{eq:bounded-transfer}
 (\Gamma'_{b,n}c)(f_0,\ldots,f_d)
  =\int_{X_n(N_g)}
       c(\gamma_n(f_0,x),\ldots,
         \gamma_n(f_d,x))\,d\mu^n(x).
\end{equation}
The cocycle identity, right homogeneity, and change of variables show that
this is a cochain map of norm at most one. It induces maps
$$
\Gamma_{b,n}^d:\Hb^d(\widehat B_n(N_g))\longrightarrow\Hb^d(\mathcal H),\qquad
\overline\Gamma_{b,n}^{\,d}:\rHb^d(\widehat B_n(N_g))\longrightarrow\rHb^d(\mathcal H).
$$
We use the same notation for their restrictions to density-preserving
diffeomorphisms. Integrating over a product
$X_{n_1}(N_g)\times\cdots\times X_{n_k}(N_g)$ gives the mixed transfer
from the corresponding product of braid groups. This is the product
construction of \cite{Brandenbursky2026}.
For density-preserving diffeomorphisms, the finite-image lemma below also
defines the transfer on
ordinary cochains and makes it commute with comparison. For homeomorphisms
we assert exactness only when a suitable primitive is integrable.

\emph{External products.}
For groups $K_1,\ldots,K_k$ and classes
$u_j\in\Hb^{d_j}(K_j)$, let
$p_j:K_1\times\cdots\times K_k\to K_j$ be the coordinate
projections. Define the external product by
$$
u_1\times\cdots\times u_k=p_1^*u_1\smile\cdots\smile p_k^*u_k
\in\Hb^{d_1+\cdots+d_k}(K_1\times\cdots\times K_k).
$$
Cup and external products descend to reduced bounded cohomology and satisfy
$$
\|u\smile v\|\le\|u\|\|v\|,\qquad\|u\times v\|\le\|u\|\|v\|.
$$
Comparison preserves cup products; hence a cup product with an exact
bounded factor is exact.

Let $\F=\langle a,b\rangle$. We use
\begin{equation}\label{eq:free-dimensions}
\dim\rHb^2(\F)=\dim\rHb^3(\F)=\infty;
\end{equation}
see \cite{Brooks1981} for degree two and
\cite[Theorem~1.2 and Sections~6--7]{Farre2018} for degree three;
see also \cite{Soma1997} for the bounded volume classes used in the latter
construction.
Since $H^q(\F)=0$ for $q\ge2$, these groups are reduced exact bounded
cohomology groups as well.

Fix $0\ne\theta\in\rHb^2(\F)$.
The degree-two duality with $\ell^1$-homology supplies a cycle
pairing nontrivially with $\theta$; see
\cite[Theorem~6.5 and Corollary~6.7]{Frigerio}.
Thus the detecting-cycle hypothesis in the cited external-product
result is satisfied.
For $m\ge0$ and $q\in\{2,3\}$, put
$$
E_m^q(v)=\theta^{\times m}\times v,\qquad
T_m^q(v)=\Sym_{m+1}(\theta^{\times m}\times v),
$$
where $\Sym_{m+1}$ is the sum of pullbacks by all permutations of the
$m+1$ factors. Here $\theta^{\times0}=1\in\rHb^0(\{1\})$.
These are maps from $\rHb^q(\F)$ to
$\rHb^{2m+q}(\F^{m+1})$.

\begin{proposition}[{\cite[Lemma~8 and Proposition~9]{Brandenbursky2026}}]
\label{prop:external-products}
The maps $E_m^q$ and $T_m^q$ are injective and admit bounded linear left
inverses. More generally, for any group $K$, external product with
$\theta^{\times m}$ gives an injective map
$$
\rHb^q(K)\longrightarrow\rHb^{2m+q}(\F^m\times K),\qquad v\longmapsto\theta^{\times m}\times v,
$$
with a bounded linear left inverse.
\end{proposition}

\section{The Gambaudo--Ghys construction for density-preserving diffeomorphisms}
\label{sec:smooth-gg}

We first prove that a fixed smooth density-preserving
diffeomorphism produces only finitely many trajectory braids,
outside a null set. The proof follows the argument of
\cite[Section~2.B.2]{Entropy}. We recall the homological
definition of flux and deduce the required fragmentation
from the corresponding result for surfaces with boundary.

\emph{Density flux and fragmentation.}
For a smooth vector field $X$ preserving $\omega$, let
$\mathcal F_\omega(X)\in H_1(N_g;\R)$ be determined by
$$
 \left\langle[\alpha],\mathcal F_\omega(X)\right\rangle
 =\int_{N_g}\alpha(X)\,d\mu
$$
for every smooth closed real one-form $\alpha$.
The integral vanishes for exact $\alpha$, so it depends
only on $[\alpha]$.

For a smooth $\omega$-preserving isotopy $f_t$ starting
at the identity, with generating vector field
$X_t=\dot f_t\circ f_t^{-1}$, define
$$
 \widetilde{\operatorname{Flux}}([f_t])
 =\int_0^1\mathcal F_\omega(X_t)\,dt
 \in H_1(N_g;\R).
$$
Here $[f_t]$ denotes the homotopy class relative to the
endpoints. This is Fathi's mass-flow homomorphism restricted
to smooth density-preserving isotopies; see
\cite[Section~5]{Fathi1980}.
Let $\Gamma_\omega$ be the subgroup consisting of the
fluxes of smooth loops based at the identity.
Passing to the quotient gives a homomorphism
$$
 \operatorname{Flux}:\Diff_0(N_g,\omega)
 \longrightarrow H_1(N_g;\R)/\Gamma_\omega.
$$

For $g>1$, write $N_g$ as a sphere with $g$ holes,
with a M\"obius band attached to each boundary circle.
Let $c_1,\ldots,c_{g-1}$ be the boundaries of $g-1$
of these bands, with chosen orientations.
These curves are pairwise disjoint, and their classes
form a basis of $H_1(N_g;\R)$.

Each curve $c_i$ is two-sided: its collars in the
M\"obius band and in the complement together form
an annular neighborhood.
Choose these annular neighborhoods pairwise disjoint
and fix an orientation on each of them.
On the $i$-th annulus, choose area coordinates
$(\theta,r)\in(\R/\Z)\times I_i$
in which $\omega$ is represented by $d\theta\wedge dr$,
with the angular direction agreeing with the orientation
of $c_i$.
Choose $b_i\in C_c^\infty(I_i)$ with
$\int_{I_i}b_i(r)\,dr=1$, and extend
$$
 Y_i=b_i(r)\partial_\theta
$$
by zero outside the annulus.
Then $Y_i$ preserves $\omega$ and
$$
 e_i:=\mathcal F_\omega(Y_i)=[c_i].
$$
Let $\phi_i^s$ denote its time-$s$ flow.
These flows commute because their supports are disjoint;
compare the annular pushes in
\cite[Lemma~2.3]{Entropy}.

\begin{lemma}\label{lem:density-fragmentation}
Let $g\ge1$. Every $h\in\ker\operatorname{Flux}$ is a finite
product $h=h_1\cdots h_s$, where each $h_i$ is a Hamiltonian
diffeomorphism compactly supported in an embedded disc in
$N_g$, after choosing an orientation on that disc.
\end{lemma}

\begin{proof}
Choose a smooth density-preserving isotopy $h_t$ from
$\id$ to $h$ with zero total flux.
Such an isotopy exists by concatenating an initial
isotopy with a loop canceling its flux.
Choose sufficiently small closed discs
$D\subset\operatorname{int}D'\subset N_g$.
By continuity and compactness of $[0,1]$, we can choose
a subdivision
$0=t_0<t_1<\cdots<t_\ell=1$
and oriented coordinate discs $U_1,\ldots,U_\ell\subset N_g$
such that
$h_t^{-1}(D')\subset U_j\quad(t\in[t_{j-1},t_j],\ 1\le j\le\ell)$.
Reparametrize each segment smoothly so that the isotopy
is stationary near its endpoints. Apply the Hamiltonian isotopy extension
construction of
\cite[proof of Theorem~3.3.1, Step~1]{McDuffSalamon2017}
in these coordinate discs. This gives a density-preserving
isotopy $k_t$, starting at the identity, which agrees
with $h_t^{-1}$ on a neighborhood of $D$.
On each subinterval, the local Hamiltonian is cut off
inside the corresponding coordinate disc.
Consequently, $k_1$ is a finite product of Hamiltonian
diffeomorphisms supported in discs, and the isotopy
$k_t$ has zero flux on $N_g$.

The isotopy $r_t=h_tk_t$ fixes a neighborhood of $D$
throughout. Put
$$
 N_{g,1}=N_g\setminus\operatorname{int}D,
 \qquad r=r_1.
$$
Let $v\in H_1(N_{g,1};\R)$ be the homological flux of
the restricted isotopy $r_t$.
Its image in $H_1(N_g;\R)$ is the flux of $r_t$ on
the closed surface, which is zero.
Since
$$
 \ker\bigl(H_1(N_{g,1};\R)\longrightarrow
                 H_1(N_g;\R)\bigr)
 =\R[\partial D],
$$
we have $v=\lambda[\partial D]$ for some $\lambda\in\R$.

Choose an annular push $a$, compactly supported in a
collar outside $D$, whose flux on $N_{g,1}$ is
$-\lambda[\partial D]$.
Viewed on $N_g$, this map is Hamiltonian in a slightly
larger disc containing $D$.
The map $ar$ belongs to the boundary-flux kernel and
is therefore a finite product of disc-supported
Hamiltonian diffeomorphisms by
\cite[Lemma~5.7]{KimMaruyamaFlux2025}.
Since
$$
 h=a^{-1}(ar)k_1^{-1},
$$
this gives the required factorization of $h$.
\end{proof}

\begin{lemma}\label{lem:autonomous-flux-correction}
Let $g>1$. For every $f\in\Diff_0(N_g,\omega)$, there is
an autonomous density-preserving diffeomorphism $u$,
supported in the union of the annuli above, such that
$uf\in\ker\operatorname{Flux}$.
\end{lemma}

\begin{proof}
Choose real numbers $s_1,\ldots,s_{g-1}$ such that
$\sum_i s_i e_i$ represents $-\operatorname{Flux}(f)$,
and set
$$
 u=\phi_1^{s_1}\circ\cdots\circ\phi_{g-1}^{s_{g-1}}.
$$
Since the supports are disjoint, $u$ is the time-one
map of the autonomous vector field $\sum_i s_iY_i$.
Its defining isotopy has flux $\sum_i s_i e_i$.
Hence $\operatorname{Flux}(u)=-\operatorname{Flux}(f)$,
which gives $\operatorname{Flux}(uf)=0$.
\end{proof}

\emph{Finitely many braids.}
A measurable function has \emph{essentially finite image}
if its image on some full-measure subset of its domain
is finite.

\begin{lemma}\label{lem:finite-braids}
For every $g\ge1$, $n\ge2$, and
$f\in\Diff_0(N_g,\omega)$, the map
$$
 X_n(N_g)\longrightarrow\widehat B_n(N_g),\qquad
 x\longmapsto\gamma_n(f,x)
$$
has essentially finite image.
\end{lemma}

\begin{proof}
For $g=1$, we have $H_1(N_1;\R)=0$, so
Lemma~\ref{lem:density-fragmentation} expresses $f$
as a finite product of Hamiltonian diffeomorphisms
supported in oriented discs.

For $g>1$, Lemmas~\ref{lem:autonomous-flux-correction}
and \ref{lem:density-fragmentation} give
$f=u^{-1}f_1\cdots f_s$,
where $u$ is a product of commuting annular pushes and
each $f_i$ is a Hamiltonian diffeomorphism supported
in an oriented disc.

Every compactly supported Hamiltonian disc map is a finite
product of Morse-type autonomous maps, by adapting the
generation argument in \cite[proof of Lemma~2.3]{Entropy},
using simplicity of the Calabi kernel in the disc case.
Here ``Morse-type'' has the meaning of
\cite[Definition~3.1]{BrandenburskyKedra2013};
in particular, the Hamiltonian is allowed to vanish
on a boundary collar. For a Morse-type autonomous disc map, the 
argument of \cite[proof of Lemma~2.1]{Entropy}, which
uses \cite[Section~2.C]{Brandenbursky2015}, gives
essentially finite image.
For $g>1$, the same argument treats the annular factors
of $u$.
Both arguments are local to the oriented supporting discs
and annuli, where the density is represented by an area form.
The closing paths are those used in that argument,
so the conclusion concerns finitely many braid
values. 

Finally, suppose the cocycles of $v$ and $w$ take
values almost everywhere in finite sets $K_v$ and $K_w$.
Equation~\eqref{eq:cocycle} implies that the cocycle
of $vw$ takes values almost everywhere in $K_vK_w$.
The cocycle of $v^{-1}$ takes values in $K_v^{-1}$,
after a change of variables.
Since all maps preserve $\mu$, the exceptional sets
remain null.
This is \cite[Proposition~2.2]{Entropy}, and it
completes the proof.
\end{proof}

\begin{proposition}\label{prop:smooth-transfer}
For $g\ge1$ and $n\ge2$, formula~\eqref{eq:gg} defines
a linear map
$$
\mathcal G_n^\omega:\Q(\widehat B_n(N_g))\longrightarrow\Q(\Diff_0(N_g,\omega)).
$$
Moreover, the smooth configuration-space transfers
preserve exact bounded cohomology in every degree,
as do their mixed product versions.
\end{proposition}

\begin{proof}
Lemma~\ref{lem:finite-braids} makes the integral in
\eqref{eq:gg} finite for every homogeneous quasimorphism.
Equation~\eqref{eq:cocycle} gives
$D(\Psi_{n,q}^{\omega})\le D(q)$. Homogenization is linear in $q$,
which proves the first assertion.
For an arbitrary ordinary homogeneous cochain $c$
and a fixed tuple $(f_0,\ldots,f_d)$ of elements of
$\Diff_0(N_g,\omega)$, the tuple of cocycle values in
\eqref{eq:bounded-transfer} has essentially finite image.
Thus the integral is finite without assuming that $c$
is bounded.
Integration commutes with the differential, and the
ordinary and bounded transfers agree on bounded cochains.
Consequently the transfer commutes with comparison.
The same argument applies to a finite product of
configuration spaces.
\end{proof}

\emph{Injectivity.}
Proposition~\ref{prop:smooth-transfer} defines the smooth
Gambaudo--Ghys map.
We will use the following naturality property.
Suppose that $h$ is a smooth diffeomorphism isotopic
to the identity and $h_*\nu=\mu$, where $\nu$ is
another smooth positive probability density, also
denoted by its associated measure.
For every $f\in\Diff_0(N_g,\omega)$, the conjugate
$h^{-1}fh$ belongs to $\Diff_0(N_g,\nu)$, and
\begin{equation}\label{eq:GG-density-naturality}
 \mathcal G_n^\omega(q)(f)
 =\mathcal G_n^\nu(q)(h^{-1}fh).
\end{equation}
Indeed, choose a smooth isotopy from $\id$ to $h^{-1}$.
For every $k\ge1$, the trajectory cocycle relation gives
$$
\gamma_n(h^{-1}f^kh,h^{-1}x)=\gamma_n(h^{-1},f^kx)\gamma_n(f^k,x)\gamma_n(h^{-1},x)^{-1}.
$$
The quasimorphism inequality therefore implies
$$
\begin{aligned}
 &\bigl|q(\gamma_n(h^{-1}f^kh,h^{-1}x))-q(\gamma_n(f^k,x))
 -q(\gamma_n(h^{-1},f^kx))\\
 &\qquad+q(\gamma_n(h^{-1},x))\bigr|\le 2D(q).
\end{aligned}
$$
Lemma~\ref{lem:finite-braids} applies to $f^k$ with
density $\omega$ and to $h^{-1}f^kh$ with density $\nu$.
Hence the first two terms are
integrable with respect to $\mu^n$.
Consequently, 
$$
 q(\gamma_n(h^{-1},f^kx))-q(\gamma_n(h^{-1},x))
$$
is also integrable.
It is integrable over the
$\mu^n$-preserving map $f^k$, so its integral vanishes:
$$
\int_{X_n(N_g)}\bigl(q(\gamma_n(h^{-1},f^kx))-q(\gamma_n(h^{-1},x))\bigr)\,d\mu^n(x)=0.
$$
Integrating the quasimorphism inequality and changing
variables by $h^{-1}$, we obtain
$$
\left|\Psi_{n,q}^{\nu}(h^{-1}f^kh)-\Psi_{n,q}^{\omega}(f^k)\right|\le 2D(q).
$$
Dividing by $k$ and passing to the limit proves
\eqref{eq:GG-density-naturality}.

The following argument adapts the injectivity arguments of
Ishida~\cite[proof of Theorem~1.2]{Ishida2014} and
Brandenbursky--Marcinkowski~\cite[Section~2.D, proof of
Theorem~2.5]{Entropy}.
We vary a smooth invariant density while keeping
the diffeomorphism fixed.

\begin{proposition}\label{prop:smooth-GG-injective}
For $g\ge1$ and $n\ge2$, the map
$$
\mathcal G_n^{\omega}:\Q(\widehat B_n(N_g))\longrightarrow\Q(\Diff_0(N_g,\omega))
$$
is injective. Moreover, if $\mathcal G_n^{\omega}(q)$ is a
homomorphism, then $q$ is a homomorphism.
\end{proposition}

\begin{proof}
Let $q\ne0$. Choose a braid $\beta\in B_n(N_g)$ on whose
image in $\widehat B_n(N_g)$ the value of $q$ is nonzero.
Replacing $\beta$ by a positive power, we may first assume
that it is pure. We then replace it by its square.
Each strand now represents the square of a loop in $N_g$,
and hence preserves local orientation.
The value of $q$ remains nonzero by homogeneity.

Choose pairwise disjoint small discs
$U_1,\ldots,U_n$ about the base points.
We realize $\beta^{-1}$ by concatenating smooth
$\omega$-preserving isotopies which move these discs
along its strands, following the construction in
\cite[Section~2.B, proof of Theorem~2]{Brandenbursky2015}.
The local motions take place in oriented coordinate
charts, where $\omega$ is represented by an area form.
The inverse appears because of our inverse-trajectory
convention.
Since every strand represents an orientation-preserving
loop, the motions are chosen so that each disc returns
with its initial parametrization.
Let $f\in\Diff_0(N_g,\omega)$ be the time-one map of
the resulting isotopy.
Then $f$ is the identity on each $U_i$ and
\begin{equation}\label{eq:realized-braid}
 q\bigl(\gamma_n(f,x)\bigr)=q(\beta)
 \qquad (x_i\in U_i).
\end{equation}
Here and below, $q(\beta)$ denotes evaluation on the image
of $\beta$ in the central quotient.
The equality is unaffected by the common conjugation
coming from the closing paths.
More generally, permuting the labels of the points
conjugates the braid in the full braid group and therefore
does not change its value under $q$.

For each $i$, choose a nonnegative smooth function
$\rho_i$ compactly supported in $U_i$ such that
$\int\rho_i\,d\mu=1$, and put
$d\nu_i=\rho_i\,d\mu$.
For $t=(t_1,\ldots,t_n)$ with $t_i>0$ and $\sum_i t_i=1$, define
$$
\nu_{0,t}=\sum_{i=1}^n t_i\nu_i,\qquad\nu_{\eps,t}
=\eps\mu+(1-\eps)\nu_{0,t},\qquad 0<\eps<1.
$$
The density of $\nu_{\eps,t}$ with respect to $\omega$ is
$\eps+(1-\eps)\sum_{i=1}^n t_i\rho_i$,
which is smooth and strictly positive.
Thus $\nu_{\eps,t}$ is a good probability measure.
The measure $\nu_{0,t}$ is used only in the polynomial
calculation below.
The map $f$ preserves $\mu$ and is the identity on the
supports of all the $\rho_i$.
Hence it preserves $\nu_{\eps,t}$.

For every $k\ge1$, define the homogeneous polynomial
of degree $n$
\begin{equation}\label{eq:GG-polynomial}
\begin{split}
&P_q(t_1,\ldots,t_n)(f^k)=\\
&\sum_{(i_1,\ldots,i_n)\in\{1,\ldots,n\}^{n}}
 t_{i_1}\cdots t_{i_n}
 \int_{X_n(N_g)}
 q\bigl(\gamma_n(f^k,x)\bigr)\,
 d(\nu_{i_1}\otimes\cdots\otimes\nu_{i_n})(x).
\end{split}
\end{equation}
The measures in \eqref{eq:GG-polynomial} are nonatomic, so the
diagonals have measure zero.

Every configuration in the support of
$\nu_{0,t}^{\otimes n}$ is fixed by $f$.
The cocycle identity therefore gives
$\gamma_n(f^k,x)=\gamma_n(f,x)^k$
for almost every such configuration.
By homogeneity of $q$ we obtain
\begin{equation}\label{eq:GG-polynomial-power}
P_q(t_1,\ldots,t_n)(f^k)=kP_q(t_1,\ldots,t_n)(f).
\end{equation}

The coefficient of $t_1\cdots t_n$ in
$P_q(t_1,\ldots,t_n)(f)$ is $n!q(\beta)$.
Indeed, the contributing terms correspond to the $n!$
assignments of the labeled points to distinct discs, and
\eqref{eq:realized-braid}, together with conjugation
invariance, gives the same value in each term.
Thus $P_q(t_1,\ldots,t_n)(f)$ is not the zero polynomial.
By homogeneity, it cannot vanish on the whole open simplex
$$
\bigl\{(t_1,\ldots,t_n):t_i>0,\ t_1+\cdots+t_n=1\bigr\}.
$$
Fix $t$ in this simplex such that
$P_q(t_1,\ldots,t_n)(f)\ne0$. Put
$$
B_q(f)=\operatorname*{ess\,sup}_{x\in X_n(N_g)}\left|q\bigl(\gamma_n(f,x)\bigr)\right|.
$$
Here $\operatorname*{ess\,sup}$ denotes the essential
supremum with respect to $\mu^n$.
This number is finite by Lemma~\ref{lem:finite-braids}.
For every $k\ge1$, the cocycle identity gives
$$
\gamma_n(f^k,x)=\gamma_n(f,f^{k-1}x)\cdots\gamma_n(f,fx)\gamma_n(f,x).
$$
Hence
$$
\left|q\bigl(\gamma_n(f^k,x)\bigr)-\sum_{j=0}^{k-1}q\bigl(\gamma_n(f,f^jx)\bigr)\right|
\le (k-1)D(q)
$$
for almost every $x$.
Since $f$ preserves $\mu^n$, the exceptional sets arising
from the $k$ factors remain null. We obtain
\begin{equation}\label{eq:Bq-power-bound}
 B_q(f^k)\le kB_q(f)+(k-1)D(q).
\end{equation}

Expanding the product measure gives
$\nu_{\eps,t}^{\otimes n}=(1-\eps)^n\nu_{0,t}^{\otimes n}+R_{\eps,t}$,
where $R_{\eps,t}$ is a positive measure, absolutely continuous
with respect to $\mu^n$, of total mass $1-(1-\eps)^n$.
Using \eqref{eq:GG-polynomial-power}, we have
$$
\frac1k\int_{X_n(N_g)}q\bigl(\gamma_n(f^k,x)\bigr)\,d\nu_{0,t}^{\otimes n}(x)=P_q(t_1,\ldots,t_n)(f).
$$
Moreover, by \eqref{eq:Bq-power-bound},
$$
\begin{aligned}
&\frac1k
\left|\int_{X_n(N_g)}q\bigl(\gamma_n(f^k,x)\bigr)\,dR_{\eps,t}(x)\right|\le\frac{B_q(f^k)}{k}\bigl(1-(1-\eps)^n\bigr)\\
&\le\left(B_q(f)+\left(1-\frac1k\right)D(q)\right)\bigl(1-(1-\eps)^n\bigr).
\end{aligned}
$$
Passing to the limit as $k\to\infty$, we obtain
$$
\left|\mathcal G_n^{\nu_{\eps,t}}(q)(f)-(1-\eps)^nP_q(t_1,\ldots,t_n)(f)\right|\le
\bigl(B_q(f)+D(q)\bigr)\bigl(1-(1-\eps)^n\bigr).
$$
Consequently,
$\mathcal G_n^{\nu_{\eps,t}}(q)(f)\ne0$
for all sufficiently small $\eps>0$.
Moser's theorem for densities \cite[Theorem~7]{BruverisMichorParusinskiRainer2018}
provides a smooth
diffeomorphism $h$, isotopic to the identity, such that
$h_*\nu_{\eps,t}=\mu$.
Equation~\eqref{eq:GG-density-naturality} gives
$$
\mathcal G_n^\omega(q)(hfh^{-1})=\mathcal G_n^{\nu_{\eps,t}}(q)(f)\ne0.
$$
This proves injectivity.

Now we prove injectivity modulo homomorphisms.
Suppose that $q$ is not a homomorphism.
Let $H$ be the image in $\widehat B_n(N_g)$ of the
finite-index subgroup of pure braids all of whose strands
preserve local orientation.
Since restriction to a finite-index subgroup detects whether a
homogeneous quasimorphism is a homomorphism, $q|_H$ is not a
homomorphism.
Hence there exists $\overline\beta\in[H,H]$ such that
$q(\overline\beta)\ne0$. Write $\overline\beta$ as a finite product of commutators
in $H$, and choose pure orientation-preserving braid lifts
of all the factors. Their commutators give a braid
$\beta\in B_n(N_g)$ whose image is $\overline\beta$.
Realize all these factors by density-preserving maps which
are the identity on the same discs
$U_1,\ldots,U_n$, and let $f$ be the corresponding product
of commutators.
Then $f$ belongs to the commutator subgroup of
$\Diff_0(N_g,\omega)$ and realizes $\beta$ on the
configurations used above.

Each factor preserves $\nu_{\eps,t}$,
because it preserves $\mu$ and is the identity on all
$U_i$.
The preceding argument gives
$\mathcal G_n^\omega(q)(hfh^{-1})\ne0$.
Since $hfh^{-1}$ belongs to the commutator subgroup of
$\Diff_0(N_g,\omega)$,
$\mathcal G_n^\omega(q)$ cannot be a homomorphism.
Thus the induced map modulo homomorphisms is injective.
\end{proof}

\begin{remark}\label{rem:GG-exact-degree-two}
Homomorphisms on the braid group are sent to homomorphisms by
\eqref{eq:gg}. Proposition~\ref{prop:smooth-GG-injective}
gives injectivity after quotienting both spaces by
homomorphisms. Under the identification of these quotients with exact
bounded cohomology in degree two, this is the map induced by the
trajectory transfer. Indeed,
$[\delta\mathcal G_n^{\omega}(q)]=[\delta\Psi_{n,q}^{\omega}]$
and $\delta\Psi_{n,q}^{\omega}$ is the transfer of $\delta q$.
\end{remark}

\section{Measure-preserving homeomorphisms}\label{sec:homeomorphisms}
Fix $g\ge2$ and $n\ge2$, and recall the space $\Q_{\mathrm{disc}}(\widehat B_n(N_g))$
defined in \eqref{eq:disc-space}.

\emph{Boundedness and continuity.}
We first assume that $\mu$ is induced by a smooth positive probability
density. The following is the local argument used in
\cite[Section~3, proof of Theorem~1.1]{BKMS2026}.

\begin{lemma}\label{lem:homeomorphism-braid-bound}
Let $q\in\Q_{\mathrm{disc}}(\widehat B_n(N_g))$. For every
$f\in\Homeo_0(N_g,\mu)$,
$$
B_q(f):=\operatorname*{ess\,sup}_{x\in X_n(N_g)}|q(\gamma_n(f,x))|<\infty.
$$
Moreover, there are a neighborhood $\mathcal U$ of the identity in the
$C^0$ topology and a constant $C_q$ such that $B_q(f)\le C_q$ for every
$f\in\mathcal U$.
\end{lemma}

\begin{proof}
Choose a finite cell decomposition with one open two-cell $\Delta$,
and use the closing paths fixed in Section~\ref{sec:preliminaries}.
For $f$ sufficiently close to the identity, choose a uniformly small
isotopy from the identity to $f$. The cutting and closing construction
in \cite[Section~3, equations~(3.1)--(3.5)]{BKMS2026} gives
$$
\gamma_n(f,x)=\alpha_{f,x}\beta_{f,x},\qquad \alpha_{f,x}\in K,\quad
\beta_{f,x}\in B_n(\Delta),
$$
outside a null set, where $K\subset\widehat B_n(N_g)$ is a fixed finite set independent of $f$
in a sufficiently small neighborhood. The construction joins short
trajectory segments in contractible coordinate neighborhoods and moves
the resulting local braids into $\Delta$. These neighborhoods are discs;
no orientation of the whole surface is used. The same construction
therefore applies to $N_g$.

Since $q$ vanishes on $B_n(\Delta)$, we obtain
$|q(\gamma_n(f,x))|\le\max_{\alpha\in K}|q(\alpha)|+D(q)$.
This proves the uniform bound on a symmetric neighborhood $\mathcal U$
of the identity. The subgroup generated by $\mathcal U$ is open and
closed, so it contains the identity component. Write an arbitrary $f$
as a finite product of elements of $\mathcal U$. Applying the cocycle
relation and the defect inequality to this product proves the lemma.
\end{proof}

\begin{proposition}\label{prop:homeomorphism-GG}
For a good probability measure $\mu$ on $N_g$, the Gambaudo--Ghys
construction defines an injective linear map
$$
\mathcal G_n^{\mu}:\Q_{\mathrm{disc}}(\widehat B_n(N_g))\longrightarrow \Q_{C^0}(\Homeo_0(N_g,\mu)),
$$
where the target denotes the space of $C^0$-continuous homogeneous
quasimorphisms. For a smooth density, its restriction to
$\Diff_0(N_g,\omega)$ is the map in
Proposition~\ref{prop:smooth-GG-injective}. It also induces an injective
map after quotienting the source and the target by their homomorphisms.
\end{proposition}

\begin{proof}
For a smooth density, Lemma~\ref{lem:homeomorphism-braid-bound} makes
\eqref{eq:gg} meaningful for every
$f\in\Homeo_0(N_g,\mu)$. The same cocycle calculation gives a
quasimorphism with defect at most $D(q)$, and
\eqref{eq:gg} defines its homogenization. By the uniform
bound in Lemma~\ref{lem:homeomorphism-braid-bound}, the integral and its
homogenization are bounded on a neighborhood of the identity.
A homogeneous quasimorphism which is bounded on a neighborhood of the
identity is continuous; see \cite[Theorem~1]{Shtern2001}, as used in
\cite[proof of Theorem~1.1(2)]{BKMS2026}.

Injectivity, including the assertion modulo homomorphisms, follows by
restriction to the diffeomorphism subgroup and
Proposition~\ref{prop:smooth-GG-injective}. For a general good measure,
apply Lemma~\ref{lem:good-measures}. Conjugation gives a topological
group isomorphism with a group preserving a smooth density; the induced
isomorphism of braid groups preserves the center and the family of
disc braid subgroups, and the proof follows.
\end{proof}

\emph{Infinite-dimensional families for every number of strands.}
\begin{proposition}\label{prop:disc-vanishing-abundance}
For every $g\ge2$ and $n\ge2$, the quotient
$$
\Q_{\mathrm{disc}}(\widehat B_n(N_g))\big/\bigl(\Q_{\mathrm{disc}}(\widehat B_n(N_g))
\cap\operatorname{Hom}(\widehat B_n(N_g),\R)\bigr)
$$
is infinite dimensional.
\end{proposition}

\begin{proof}
The point-pushing homomorphism identifies $\widehat B_n(N_g)$ with an
irreducible subgroup of $\Mod(N_g,n)$ which is not virtually abelian.
The Bestvina--Fujiwara construction, applied to this
subgroup, gives an infinite-dimensional space of
homogeneous quasimorphisms on $\widehat B_n(N_g)$ which vanish on
reducible elements; see \cite[Theorem~12]{BestvinaFujiwara2002} and
\cite[proof of Theorem~1.1 and Theorem~5.15]{BestvinaFujiwara2007}.

Every element of $\iota_D(B_n(D))$ is reducible,
since it preserves the essential curve $\partial D$.
Consequently, all these quasimorphisms belong to
$\Q_{\mathrm{disc}}(\widehat B_n(N_g))$, and hence
$\dim\Q_{\mathrm{disc}}(\widehat B_n(N_g))=\infty$.
Since $\widehat B_n(N_g)$ is finitely generated,
$\operatorname{Hom}(\widehat B_n(N_g),\R)$ is finite dimensional.
Therefore
$$
\dim\left(\Q_{\mathrm{disc}}(\widehat B_n(N_g))\big/\bigl(\Q_{\mathrm{disc}}(\widehat B_n(N_g))
\cap\operatorname{Hom}(\widehat B_n(N_g),\R)\bigr)\right)=\infty.
$$
\end{proof}

Together, Propositions~\ref{prop:homeomorphism-GG} and
\ref{prop:disc-vanishing-abundance} prove that, for every 
$n\ge2$, the image of the homeomorphism Gambaudo--Ghys map is
infinite dimensional, also modulo homomorphisms. In particular, its
degree-two exact bounded cohomology image is infinite dimensional.
There are no nonzero classes of zero seminorm in degree two, so the
same assertion holds in reduced exact bounded cohomology.

For completeness, the smooth construction also has infinite-dimensional
image modulo homomorphisms when $g=1$ and $n\ge3$. In this case
$\widehat B_n(N_1)$ is a finite-index subgroup of $\Mod(N_1,n)$,
and \cite[Theorems~1.1 and~5.15]{BestvinaFujiwara2007} apply.
For $n=2$, the group $B_2(N_1)$ is finite, so its space of homogeneous
quasimorphisms is zero.

\begin{remark}\label{rem:Calabi-reflection}
Assume that $\mu$ is induced by a smooth density. Let $D$ be an embedded
disc, and suppose that a homogeneous quasimorphism on
$\Homeo_0(N_g,\mu)$ restricts to $c\operatorname{Cal}_D$ on
$\Ham_c(D)$. An isotopy which moves $D$ once around a one-sided loop
gives a density-preserving map whose restriction to $D$ reverses its
orientation. If this restriction is denoted by $R$,
then $R^*\omega_D=-\omega_D$ for the corresponding area form.
If $f$ is generated by $F_t$, the map $RfR^{-1}$ is generated by
$-F_t\circ R^{-1}$, and hence
$\operatorname{Cal}_D(RfR^{-1})=-\operatorname{Cal}_D(f)$.
Conjugation invariance therefore forces $c=0$.

The same disc motion gives a useful observation for two strands. Place
the two base points on the reflection axis, and choose the final
restriction to $D$ to be that reflection. The resulting point-pushing
mapping class conjugates the disc half-twist to its inverse. Since
point-pushing identifies the braid central quotient with its image,
the image of the generator of $B_2(D)$ is conjugate to its inverse in
$\widehat B_2(N_g)$. Every homogeneous quasimorphism vanishes on this
generator. Thus, for $g\ge2$,
$\Q_{\mathrm{disc}}(\widehat B_2(N_g))=\Q(\widehat B_2(N_g))$.
\end{remark}

\section{Periodic braids and undistortion}\label{sec:undistortion}

In this section we work with $\Homeo_0(N_g)$ and $\Diff_0(N_g)$,
without requiring preservation of a measure. A smooth positive
probability density $\omega$, with associated measure $\mu$, is used
only to define the integrals. We follow the semi-bounded-cochain
and pointwise arguments of \cite[Sections~3--4]{BKMS2026}.

Recall that an element $f$ of a finitely generated group $\Gamma$ is
\emph{undistorted} if
$$
 \lim_{k\to\infty}\frac{|f^k|_S}{k}>0,
$$
where $S$ is a finite symmetric generating set of $\Gamma$.
An element of an arbitrary group
is undistorted if it is undistorted in every finitely generated subgroup
containing it.

\emph{Trajectory bounds without measure preservation.}
We use the geodesic closing paths from
Section~\ref{sec:preliminaries}.
At exceptional configurations, complete this choice
measurably by following a previously defined closing
path to a nearby configuration and then joining
corresponding points by short geodesics in pairwise
disjoint convex coordinate discs.
This adds at most one segment to each strand and
preserves the uniform bounds on the number and lengths
of segments.
The same trajectory construction therefore defines
$\gamma_n$ for arbitrary homeomorphisms, with
\eqref{eq:cocycle} holding at every configuration.
For maps preserving the fixed smooth measure, the
completed cocycle agrees almost everywhere with the
previous one, so the integrals remain unchanged.

\begin{lemma}\label{lem:ambient-braid-bounds}
For these closing paths, the following statements hold.
\begin{enumerate}[label=\textup{(\roman*)}]
\item
If $g\ge2$, $n\ge2$ and
$q\in\Q_{\mathrm{disc}}(\widehat B_n(N_g))$, then
$$
 B_q(f):=\sup_{x\in X_n(N_g)}|q(\gamma_n(f,x))|<\infty
 \qquad(f\in\Homeo_0(N_g)).
$$
For $g>2$ and $n=1$, the same bound holds for every
$q\in\Q(\pi_1(N_g))$.
\item
For $g\ge1$, $n\ge2$ and $f\in\Diff_0(N_g)$, the set
$$
 \{\gamma_n(f,x):x\in X_n(N_g)\}
 \subset\widehat B_n(N_g)
$$
is finite. In particular, $B_q(f)<\infty$ for every
$q\in\Q(\widehat B_n(N_g))$.
\end{enumerate}
\end{lemma}

\begin{proof}
We use the argument in
\cite[proof of Lemma~3.1]{BKMS2026}.
Although stated there for orientable surfaces, the same
argument applies to $N_g$: braids whose strands consist
of a uniformly bounded number of geodesic segments of
uniformly bounded length have only finitely many images
in $\widehat B_n(N_g)$.

For part~\textup{(i)}, apply the proof of
Lemma~\ref{lem:homeomorphism-braid-bound} to a sufficiently
$C^0$-small homeomorphism $f$.
The local cutting and closing construction does not use
measure preservation and gives the required bound whenever
the original closing paths at $x$ and $fx$ are defined.

For any other configuration $x$, choose a nearby
configuration $y$ for which both endpoint closing paths
are defined. Compare the trajectories using short
connectors in pairwise disjoint coordinate discs at
the two endpoints, chosen so that the image under $f$
of the first connector is homotopic to the second.
The cocycle values at $x$ and $y$ then differ by
multiplication on the left and right by loops formed
from closing paths and these connectors.
The preceding finite-set argument applies to these loops,
so the quasimorphism inequality extends the bound to
every configuration.

Every element of $\Homeo_0(N_g)$ is a finite product
of such small homeomorphisms.
The pointwise cocycle identity gives
\begin{equation}\label{eq:braid-bound-product}
 B_q(uv)\le B_q(u)+B_q(v)+D(q),
\end{equation}
which proves the assertion for arbitrary $f$.
For $n=1$, the disc braid group is trivial, so the same
argument requires no vanishing condition on $q$.

For part~\textup{(ii)}, use the same reduction to small
factors, choosing them sufficiently $C^1$-close to the
identity. Such factors admit isotopies with short geodesic
trajectories. Indeed, for a sufficiently $C^1$-small factor $v$, the
short-geodesic interpolation
$v_t(x)=\exp_x\bigl(t\exp_x^{-1}(v(x))\bigr),\thinspace 0\le t\le1$,
remains $C^1$-close to the identity uniformly in $t$.
It therefore consists of diffeomorphisms and joins
$\id$ to $v$, with each point tracing a short
geodesic segment. The argument of
\cite[proof of Lemma~3.1]{BKMS2026} shows that
the trajectory cocycle of each factor has finite image
in $\widehat B_n(N_g)$.
Thus no vanishing condition on $q$ is needed.
The cocycle identity completes the proof as in
part~\textup{(i)}.
\end{proof}

The supremum in this lemma is an ordinary supremum. We use this
meaning of $B_q$ throughout the present section.
The completion on the exceptional set is needed for the ambient
homeomorphism group, but it does not affect any of the
measure-preserving constructions in the preceding sections.

\emph{Semi-bounded cocycles.}
For a function $\psi:K\to\R$, write its differential in
inhomogeneous notation as
$(\delta\psi)(u,v)=\psi(u)+\psi(v)-\psi(uv)$.
The differential is called \emph{semi-bounded} if, for every fixed $u\in K$,
$\sup_{v\in K}|(\delta\psi)(u,v)|<\infty$.
The bound may depend on $u$.
Let $\mathcal T$ be $\Homeo_0(N_g)$, with $g\ge2$ and
$q\in\Q_{\mathrm{disc}}(\widehat B_n(N_g))$, or let
$\mathcal T=\Diff_0(N_g)$, with $g\ge1$ and arbitrary
$q\in\Q(\widehat B_n(N_g))$.
For a fixed configuration $x$, put
$$
 \Psi_{n,q,x}(f)=q(\gamma_n(f,x)),
 \qquad f\in\mathcal T.
$$
We also use the integral $\Psi_{n,q}^{\omega}$ from
\eqref{eq:gg}, with the fixed reference measure $\mu^n$.
Both functions vanish at the identity.

\begin{proposition}\label{prop:ambient-semibounded}
The differentials of the functions $\Psi_{n,q,x}$ and $\Psi_{n,q}^{\omega}$ are
semi-bounded. More precisely, for every $f,h\in\mathcal T$,
$$
 |(\delta\Psi_{n,q,x})(f,h)|\le2B_q(f)+D(q),\qquad
 |(\delta\Psi_{n,q}^{\omega})(f,h)|\le2B_q(f)+D(q).
$$
For $g>2$ and $n=1$, the same assertions hold on
$\Homeo_0(N_g)$ for every $q\in\Q(\pi_1(N_g))$.
\end{proposition}

\begin{proof}
Lemma~\ref{lem:ambient-braid-bounds} gives pointwise boundedness
and hence integrability with respect to the probability measure
$\mu^n$. The cocycle identity and the quasimorphism inequality give
$$
\begin{aligned}
&\left|q(\gamma_n(fh,y))-q(\gamma_n(h,y))-q(\gamma_n(f,y))\right|\\
&\le\left|q(\gamma_n(f,hy))-q(\gamma_n(f,y))\right|+D(q)\le2B_q(f)+D(q).
\end{aligned}
$$
Evaluating at $y=x$ proves the first estimate. Integrating with
respect to $\mu^n$ proves the second. There is no change of
variables, so $h$ need not preserve $\mu$. This is the calculation
of \cite[proof of Theorem~1.1(1) and Proposition~4.1]{BKMS2026}.
\end{proof}

We recall the group-theoretic consequence of semi-boundedness;
see \cite[Theorem~1.6]{BKMS2026}.
If $\psi(\id)=0$ and $\delta\psi$ is semi-bounded, then
for every finite symmetric generating set $S$ of a subgroup
$\Gamma\le K$ there is $C_S>0$ such that
$$
 |\psi(w)|\le C_S|w|_S\qquad(w\in\Gamma).
$$
Indeed, for each $s\in S$ the differences
$\psi(sv)-\psi(v)$ are bounded uniformly in $v$ by
$|\psi(s)|+\sup_v|\delta\psi(s,v)|$.
Write $f=s_1\cdots s_m$, where $s_i\in S$ and $m=|f|_S$. Then
$$
\begin{aligned}
|\psi(f)|
&\le
\sum_{i=1}^{m}
\left|\psi(s_i\cdots s_m)-\psi(s_{i+1}\cdots s_m)\right|\le\sum_{i=1}^{m}
\left(|\psi(s_i)|+\sup_{v\in\Gamma}|\delta\psi(s_i,v)|\right)\\
&\le m\max_{s\in S}\left(|\psi(s)|+\sup_{v\in\Gamma}|\delta\psi(s,v)|\right)\le C_S|f|_S.
\end{aligned}
$$
Here the empty product is the identity, and $C_S>0$ is chosen at least
as large as the displayed maximum.
Consequently, if
$\limsup_{k\to\infty}\frac{|\psi(f^k)|}{k}>0$,
then the stable word length of $f$ is positive in every such
subgroup containing it.

\begin{proof}[Proof of Theorem~\ref{thm:ambient-cochains}]
The existence and semi-boundedness assertions are
Proposition~\ref{prop:ambient-semibounded}.
Apply the preceding observation to $\psi=\Psi_{n,q}^{\omega}$
to obtain the undistortion assertion.
\end{proof}

\emph{Periodic configurations.}
Let $f\in\Homeo_0(N_g)$, and suppose that
$x=(x_1,\ldots,x_n)\in X_n(N_g)$ is fixed by $f^r$, where
$r\ge1$. Closing its trajectory using the fixed path
$\alpha_x$ and the inverse-trajectory convention gives
$\beta_{f,x,r}\in\widehat B_n(N_g)$.
Changing this path conjugates the braid, and changing the isotopy
does not change its image in the central quotient.
Since $f^r(x)=x$, the cocycle identity gives
\begin{equation}\label{eq:periodic-braid-powers}
 \gamma_n(f^{rk},x)=\gamma_n(f^r,x)^k=(\beta_{f,x,r})^k\qquad(k\ge1).
\end{equation}
In particular,
$|\Psi_{n,q,x}(f^{rk})|=k|q(\beta_{f,x,r})|$.

\begin{proof}[Proof of Theorem~\ref{thm:undistortion-main}]
Use $\mathcal T=\Homeo_0(N_g)$ in part~\textup{(i)} and
$\mathcal T=\Diff_0(N_g)$ in part~\textup{(ii)}.
Let $\Gamma\le\mathcal T$ be any finitely generated subgroup
containing $f$, and let $S$ be a finite symmetric generating set.
By Lemma~\ref{lem:ambient-braid-bounds}, $B_q(s)<\infty$ for every $s\in S$.
Since $S$ is finite and $D(q)<\infty$, we may choose
a constant $C_S>0$ such that
$B_q(s)+D(q)\le C_S$ for every $s\in S$.
Repeated application of \eqref{eq:braid-bound-product} to a shortest
word representing $f^{rk}$, together with
\eqref{eq:periodic-braid-powers}, gives
$$
k|q(\beta_{f,x,r})|\le B_q(f^{rk})\le C_S|f^{rk}|_S\qquad(k\ge1).
$$
Therefore
$$
\lim_{k\to\infty}\frac{|f^{rk}|_S}{rk}\ge\frac{|q(\beta_{f,x,r})|}{rC_S}>0.
$$
The stable word length of $f$ is the same limit along all positive
integers, by subadditivity. Hence $f$ is undistorted in
$\mathcal T$. Equivalently, one can apply the semi-bounded criterion
to the pointwise cochain $\Psi_{n,q,x}$.
\end{proof}

\begin{remark}
Only one periodic configuration is needed, and neither an invariant
measure nor nonzero asymptotic growth of the integrated cochain is
required. The pointwise cochain detects the return braid directly.
This is the mechanism of
\cite[Corollary~1.8 and Proposition~4.2]{BKMS2026}.
\end{remark}

\begin{corollary}\label{cor:noncontractible-periodic-point}
Let $g>2$ and $f\in\Homeo_0(N_g)$.
If $f$ has a periodic point with noncontractible return trajectory,
then $f$ is undistorted in $\Homeo_0(N_g)$.
In particular, every periodic point of a distorted element has a
contractible return trajectory.
\end{corollary}

\begin{proof}
Let $x$ have period $r$, choose $x$ as the base point, and let
$\beta\in\pi_1(N_g,x)$ be its return trajectory under $f^r$.
The group $\pi_1(N_g)$ is torsion-free and hyperbolic, and no
nontrivial element has a nonzero power conjugate to its inverse.
If $\beta\ne1$, the counting-quasimorphism construction of
\cite[Theorem~A$'$]{CalegariFujiwara2010} gives a homogeneous
quasimorphism $q$ with $q(\beta)\ne0$.
Apply the one-point assertion of
Proposition~\ref{prop:ambient-semibounded}. Since
$\Psi_{1,q,x}(f^{rk})=kq(\beta)$, the semi-bounded criterion
proves undistortion. No invariant measure is required.
\end{proof}

\emph{Infinite-order return braids for diffeomorphisms.}
For smooth maps, the finiteness assertion in
Lemma~\ref{lem:ambient-braid-bounds}\textup{(ii)} gives a criterion
which only involves the order of the braid.

\begin{theorem}\label{thm:smooth-undistortion}
Let $g\ge1$ and $f\in\Diff_0(N_g)$.
Suppose that $x\in X_n(N_g)$ is fixed by $f^r$ for some
$n\ge2$ and $r\ge1$.
If $\beta_{f,x,r}$ has infinite order in $\widehat B_n(N_g)$,
then $f$ is undistorted in $\Diff_0(N_g)$.
\end{theorem}

\begin{proof}
Using point-pushing, identify $\widehat B_n(N_g)$
with a subgroup of $\Mod(N_g,n)$, as recalled in
Section~\ref{sec:preliminaries}.
Fix a word norm $|\cdot|$ on $\Mod(N_g,n)$ and put
$$
L(v)=\sup_{y\in X_n(N_g)}|\gamma_n(v,y)|,\qquad v\in\Diff_0(N_g),
$$
where the braids are regarded as elements
of $\Mod(N_g,n)$.
By Lemma~\ref{lem:ambient-braid-bounds}\textup{(ii)},
$L(v)<\infty$ for every $v$.
The cocycle identity gives
$$
L(uv)\le\sup_{y\in X_n(N_g)}\bigl(|\gamma_n(u,vy)|+|\gamma_n(v,y)|\bigr)\le L(u)+L(v).
$$
Every infinite-order element of $\Mod(N_g,n)$ has
linear word growth; this follows from
\cite[Theorem~1.1]{KatayamaKuno2024} and
\cite[Theorem~1.2]{FarbLubotzkyMinsky2001}.
Since point-pushing is injective and $\beta_{f,x,r}$
has infinite order, there is $c>0$ such that
$|\beta_{f,x,r}^{\,k}|\ge ck$ for $k\ge 1$.

Let $\Gamma\le\Diff_0(N_g)$ be a finitely generated
subgroup containing $f$, and let $S$ be a finite
symmetric generating set.
Since $S$ is finite and $L(s)<\infty$ for every $s\in S$,
choose $C_S>0$ such that $L(s)\le C_S$ for all $s\in S$.
Equation~\eqref{eq:periodic-braid-powers} and
subadditivity of $L$ imply
$$
ck\le |\beta_{f,x,r}^{\,k}|=|\gamma_n(f^{rk},x)|\le L(f^{rk})\le C_S|f^{rk}|_S\qquad(k\ge1).
$$
Consequently,
$$
\lim_{k\to\infty}\frac{|f^{rk}|_S}{rk}\ge\frac{c}{rC_S}>0.
$$
Thus $f$ has positive stable word length in every
finitely generated subgroup containing it, and is
therefore undistorted.
\end{proof}

\begin{corollary}\label{cor:smooth-distortion-braids}
Let $g\ge1$ and let $f\in\Diff_0(N_g)$ be distorted.
For every $n\ge2$, the braid of every $n$-point periodic
configuration has finite order in $\widehat B_n(N_g)$.
If $g>2$, every such braid is trivial in $B_n(N_g)$.
\end{corollary}

\begin{proof}
The first assertion follows from
Theorem~\ref{thm:smooth-undistortion}.
For $g>2$, the group $B_n(N_g)$ is torsion-free and has trivial
center; see \cite[Propositions~1.5 and~1.6]{ParisRolfsen1999}.
Thus $\widehat B_n(N_g)=B_n(N_g)$, and the second assertion follows.
\end{proof}

\begin{remark}\label{rem:Franks-Handel-comparison}
The one-point conclusion of
Corollary~\ref{cor:noncontractible-periodic-point}
also follows from the linear-displacement argument in
\cite[Lemmas~5.4 and~6.1]{FranksHandel2006}.
This argument applies to homeomorphisms, without any
invariant-measure assumption, and works without change
on $N_g$ for $g>2$.

We are not aware of an earlier statement of
Theorem~\ref{thm:smooth-undistortion} in this generality.
For every $g\ge1$, it gives an undistortion criterion
in $\Diff_0(N_g)$ whenever a periodic configuration
has a braid of infinite order in $\widehat B_n(N_g)$,
without any invariant-measure assumption.
It implies the finite-order conclusion
of Corollary~\ref{cor:smooth-distortion-braids}.
It seems plausible that an alternative proof could
be obtained by adapting the 
arguments of Franks--Handel~\cite{FranksHandel2006}
to nonorientable surfaces.
We do not carry out this adaptation here.
Our proof uses smooth braid finiteness and linear
word growth in the punctured mapping class group.

To the best of our knowledge, the undistortion criteria
proved here for elements of $\Homeo_0(N_g)$ involving
periodic configurations with at least two points are
also new for nonorientable surfaces.
They follow from the semi-bounded cochains constructed
in this paper and require neither differentiability
nor preservation of a measure.
\end{remark}

\section{Braid groups of the projective plane and the Klein bottle}
\label{sec:low-genus}

We prove two lemmas about surface braid groups which will be used in
Section~\ref{sec:bounded-applications}. A bar over a braid denotes its
image in the central quotient.

\begin{lemma}[The projective plane]\label{lem:rp2-source}
Consider $\widehat B_3(N_1)$. Choose a disc containing the three base
points, and let $\sigma_1,\sigma_2$ be its standard Artin generators.
The homomorphism
$$
 \jmath_1:\F\longrightarrow \widehat B_3(N_1),\qquad
 \jmath_1(a)=\overline{\sigma_1^2},\quad
 \jmath_1(b)=\overline{\sigma_2^2},
$$
is injective. For $q=2,3$, the restriction map
$\jmath_1^*:\rHb^q(\widehat B_3(N_1))\longrightarrow\rHb^q(\F)$
has infinite-dimensional image. Moreover,
$\rHb^q(\widehat B_3(N_1))=\rEH^q(\widehat B_3(N_1))$ for every $q\ge2$.
\end{lemma}

\begin{proof}
The group $P_3(N_1)$ is isomorphic to $F_2\rtimes Q_8$,
where $Q_8$ is the quaternion group of order eight,
and the center of $B_3(N_1)$ is generated by the full
twist $\Delta_3^2$, of order two; see
\cite[Section~2.3]{GoncalvesGuaschi2004} and
\cite[Section~3]{GuaschiJuanPineda2026}.
Thus $\widehat B_3(N_1)$ is a non-elementary virtually
free group. Restriction to a free subgroup of finite
index is injective in ordinary real cohomology, by transfer. Consequently,
\begin{equation}\label{eq:rp2-ordinary}
 H^q(\widehat B_3(N_1);\R)=0\qquad(q\ge2).
\end{equation}
In particular, all bounded classes in these degrees
are exact.

We next prove injectivity of $\jmath_1$.
Under point-pushing, the images of $\sigma_1^2$ and
$\sigma_2^2$ correspond to Dehn twists about the standard
curves enclosing the pairs
$\{z_1,z_2\}$ and $\{z_2,z_3\}$, respectively.
These are generic two-sided curves in the three-punctured
projective plane: each bounds a twice-punctured disc,
and its complement is a once-punctured M\"obius band.
Their geometric intersection number is two.
By \cite[Theorem~13.2]{Stukow2017}, the corresponding
Dehn twists freely generate a group of rank two.
It follows that $\jmath_1$ is injective.

Since $\widehat B_3(N_1)$ is hyperbolic and
$\jmath_1(\F)$ is non-elementary,
\cite[Theorem~1.3]{HullKapovich2018}
gives a free subgroup $F'\cong F_2$ contained in
$\jmath_1(\F)$ such that $F'\times K$ is hyperbolically
embedded in $\widehat B_3(N_1)$ for some finite
subgroup $K$.
The exact extension theorem
\cite[Corollary~2 and Proposition~5.1]
{FrigerioPozzettiSisto2015}
gives a surjection
$$
 EH_b^q(\widehat B_3(N_1))
 \longrightarrow EH_b^q(F'\times K),
 \qquad q=2,3.
$$
Restriction from $F'\times K$ to $F'$ is an isometric
isomorphism in bounded cohomology, since the finite
factor $K$ is amenable.
Moreover, $H^q(F';\R)=0$ for $q\ge2$.
Passing to reduced cohomology therefore gives a
surjection
$$
 \rHb^q(\widehat B_3(N_1))
 \longrightarrow\rHb^q(F'),
 \qquad q=2,3.
$$
Since $F'\subset\jmath_1(\F)$, this restriction factors
through $\jmath_1^*$.
By \eqref{eq:free-dimensions}, $\rHb^q(F')$ is infinite
dimensional for $q=2,3$.
Hence the image of $\jmath_1^*$ is infinite dimensional.
\end{proof}

\begin{lemma}[The Klein bottle]\label{lem:klein-source}
The following statements hold.
\begin{enumerate}[label=(\roman*)]
\item Let $u,v$ be the standard generators of
$\pi_1(N_2\setminus\{z_1\},z_2)$, where $u$ preserves local
orientation and $v$ reverses it. Regard them as pure braids
which move $z_2$ while fixing $z_1$.
The homomorphism
\begin{equation}\label{eq:klein-free-lifts}
 \jmath_2:\F\longrightarrow\widehat B_2(N_2),\qquad
 \jmath_2(a)=\overline{v^2},\quad
 \jmath_2(b)=\overline{vuv^{-1}u^{-1}},
\end{equation}
is injective. Its image is represented by braids with one
stationary strand and one strand tracing an orientation-preserving
loop. For $q=2,3$, restriction gives a map
$$
 \jmath_2^*:\rEH^q(\widehat B_2(N_2))
               \longrightarrow\rHb^q(\F)
$$
with infinite-dimensional image.
Every homogeneous quasimorphism on $\widehat B_2(N_2)$
vanishes on the disc braid subgroup $B_2(D)$.

\item Let $\jmath_D:\F\to\widehat B_3(N_2)$ send $a,b$ to the
images of $\sigma_1^2,\sigma_2^2$ in a disc containing the three
base points. There is a class
$\alpha_2\in\rEH^2(\widehat B_3(N_2))$ such that
$$
 \eta_2:=\jmath_D^*\alpha_2\ne0
 \quad\text{in }\rHb^2(\F).
$$
\end{enumerate}
\end{lemma}

\begin{proof}
We first prove~\textup{(i)}. The moving-point subgroup $F(u,v)$
embeds in $P_2(N_2)$; see
\cite[Proposition~3.1]{GoncalvesGuaschiLaass2020}.
Any element of $F(u,v)\cap Z(B_2(N_2))$ belongs to the
center of $F(u,v)$, which is trivial.
Thus $F(u,v)$ also embeds in $\widehat B_2(N_2)$.
The words $v^2$ and $vuv^{-1}u^{-1}$ do not commute in $F(u,v)$
and hence freely generate a group of rank two.
This proves injectivity of $\jmath_2$.
Both words contain an even number of occurrences of $v^{\pm1}$,
so their trajectories preserve local orientation.

We show that $\jmath_2(\F)$ contains a pseudo-Anosov point-push.
The loop
$$
w=v^2(vuv^{-1}u^{-1})
$$
fills $N_2\setminus\{z_1\}$.
Let us explain this. After filling the puncture, use the same letters for
the images of $u,v$ in
$$
\pi_1(N_2)=\langle u,v\mid vuv^{-1}=u^{-1}\rangle.
$$
The image of $w$ is $v^2u^{-2}=u^{-2}v^2$, since $v^2$
commutes with $u$ in $\pi_1(N_2)$.

Suppose that $w$ can be represented by a loop disjoint
from an essential nonperipheral simple closed curve
in $N_2\setminus\{z_1\}$.
After filling the puncture, this curve remains noncontractible,
and cutting $N_2$ along it gives an annulus or M\"obius
bands. In the annulus, loops represent powers of $u$;
in either M\"obius band, orientation-preserving loops
represent powers of $v^2$, up to conjugacy.
Thus the image $u^{-2}v^2$ of $w$ would be conjugate
to $u^m$ or $v^{2m}$ for some $m\in\Z$. The homomorphism
$$
\pi_1(N_2)\longrightarrow\Z,\qquad u\longmapsto0,\quad v\longmapsto1
$$
takes $u^{-2}v^2$ to $2$ and every power of $u$ to $0$,
ruling out the first possibility.
If $u^{-2}v^2$ were conjugate to $v^{2m}$, the same
homomorphism would give $m=1$.
But $v^2$ is central, since
$v^2uv^{-2}=v u^{-1}v^{-1}=u$.
Conjugacy to $v^2$ would therefore imply
$u^{-2}v^2=v^2$, and hence $u^2=1$, contradicting
the infinite order of $u$.
Thus neither possibility occurs, and $w$ fills
the punctured Klein bottle.
The point-pushing criterion for the punctured Klein bottle
therefore makes $\jmath_2(ab)$ pseudo-Anosov; see
\cite[Theorem~3.5 and Remark~3.6]{Boke2025}.

The mapping class group $\Mod(N_2,2)$ acts coboundedly on its
hyperbolic curve graph, and this action is WPD by
\cite[Theorem~5.15]{BestvinaFujiwara2007}.
Since $\jmath_2(\F)$ contains a pseudo-Anosov element and is not
virtually cyclic, its action is non-elementary.
Apply \cite[Theorem~1.3]{HullKapovich2018} to this subgroup.
It gives a free subgroup $F'\cong F_2$ contained in
$\jmath_2(\F)$ and a finite subgroup $K$ such that
$F'\times K$ is hyperbolically embedded in $\Mod(N_2,2)$.
As in the proof of Lemma~\ref{lem:rp2-source}, the exact extension
theorem \cite[Corollary~2 and Proposition~5.1]
{FrigerioPozzettiSisto2015}, followed by restriction to $F'$
and reduction, gives a surjection
$$
 \rEH^q(\Mod(N_2,2))\longrightarrow\rHb^q(F'),
 \qquad q=2,3.
$$
This restriction factors through
$\rEH^q(\widehat B_2(N_2))$ and then through $\jmath_2^*$.
Since $\rHb^q(F')$ is infinite dimensional for $q=2,3$,
the image of $\jmath_2^*$ is infinite dimensional in both degrees.
Finally, the disc half-twist is conjugate to its inverse in
$\widehat B_2(N_2)$, as observed in
Remark~\ref{rem:Calabi-reflection}.
Every homogeneous quasimorphism therefore vanishes on it,
and hence on the cyclic disc braid subgroup $B_2(D)$.

We now prove~\textup{(ii)}. Consider the homomorphisms
$$
\widehat B_3(N_2)\longrightarrow\Mod(N_2,3)\longrightarrow\Mod^+(\Sigma_1,6)
$$
given by point-pushing and the orientation-preserving lift.
The first homomorphism was recalled in
Section~\ref{sec:preliminaries}; the second is injective by
\cite[Theorem~13(b)]{GuaschiJuanPineda2026}.
In the three-point disc, take the braid
$\beta=(\sigma_2^{-1}\sigma_1\sigma_2^{-3}\sigma_1^2)^6(\sigma_1\sigma_2)^3$.
Its exponent sum is zero, so
$\beta\in[B_3(D),B_3(D)]$.
The proof of \cite[Lemmas~4 and~5]{BrandenburskyNonorientable2026}
constructs a homogeneous quasimorphism on the lifted mapping class
group which is nonzero on this braid.
The same proof applies here: the two lifts of $D$ have essential
boundaries in the six-punctured torus, and the non-reversibility
and stretch-factor arguments for the two pseudo-Anosov restrictions
are unchanged. The required quasimorphism is supplied by
\cite[Theorem~4.2 and Proposition~4.4]{BBF2016}, as in that proof.

Pulling this quasimorphism back to $\widehat B_3(N_2)$ gives a
homogeneous quasimorphism whose restriction to $B_3(D)$ is not a
homomorphism, since it is nonzero on $\beta$.
The resulting exact degree-two bounded class therefore has nonzero
restriction to $B_3(D)$.
Restriction to $P_3(D)$ is injective in bounded cohomology, and
$$
P_3(D)\cong\langle\sigma_1^2,\sigma_2^2\rangle\times\Z\cong\F\times\Z.
$$
Removing the amenable central factor shows that its restriction
along $\jmath_D$ is nonzero.
There are no nonzero degree-two bounded cohomology classes of zero
seminorm; hence this restriction remains nonzero after reduction.
This gives $\alpha_2$ and $\eta_2$.
\end{proof}

\section{Proof of Theorem~\ref{thm:bounded-main}}
\label{sec:bounded-applications}
We first use the probability measure $\mu$ induced by the smooth
positive density $\omega$. Lemma~\ref{lem:good-measures} then gives
the assertions for every good measure. For $g=1,2$, let
$\mathcal T_g$ denote either $\Homeo_0(N_g,\mu)$ or
$\Diff_0(N_g,\omega)$.

\begin{proof}[Proof of Theorem~\ref{thm:bounded-main} for $N_1$]
Use the group $\widehat B_3(N_1)$ and the homomorphism $\jmath_1$ from
Lemma~\ref{lem:rp2-source}. For $q=2,3$, choose an
infinite-dimensional linear subspace $V_q\subset\rHb^q(\widehat B_3(N_1))$
on which $\jmath_1^*$ is injective. Choose a nonzero class
$\alpha_1\in V_2$ and put
$$
\eta_1=\jmath_1^*\alpha_1\in\rHb^2(\F),\qquad A_1=\overline\Gamma_{b,3}^{\,2}(\alpha_1)\in\rHb^2(\mathcal T_1).
$$
Thus $\eta_1\ne0$. For $m\ge0$ and $q=2,3$, define
$$
\Phi_{1,m,q}:V_q\longrightarrow\rHb^{2m+q}(\mathcal T_1),\qquad
\Phi_{1,m,q}(v)=A_1^{\smile m}\smile\overline\Gamma_{b,3}^{\,q}(v).
$$
Fubini identifies this cup product with the transfer of
$\alpha_1^{\times m}\times v$ over $X_3(N_1)^{m+1}$;
see \cite[Section~4.1]{Brandenbursky2026}.

Fix $m$ and put $k=m+1$. Choose pairwise disjoint discs
$D_1,\ldots,D_k\subset N_1$. In an oriented area coordinate
chart in each $D_j$, choose two overlapping round discs
$D_{j,1},D_{j,2}\Subset D_j$ such that their intersection
and both exclusive regions have positive measure.
For $s=1,2$, let $f_{j,s,\eps}$ be the endpoint of a smooth
rotation which makes one full turn on a slightly smaller
concentric disc and is the identity near $\partial D_{j,s}$.
Choose the transition collars with total measure tending
to zero as $\eps\to0$. The direction of rotation is chosen
so that, for a configuration outside the transition collars with one point
in each of
$$
 D_{j,1}\setminus D_{j,2},\qquad
 D_{j,1}\cap D_{j,2},\qquad
 D_{j,2}\setminus D_{j,1},
$$
the inverse trajectories of the two rotations represent
$\sigma_1^2,\sigma_2^2$, respectively. These are the local
maps of \cite[Lemmas~15 and~21]{Brandenbursky2026}.
Define
$$
 \rho_\eps:\F^k\longrightarrow\Diff_0(N_1,\omega),\qquad
 \rho_\eps(w_1,\ldots,w_k)
     =\prod_{j=1}^k w_j(f_{j,1,\eps},f_{j,2,\eps}).
$$
Here each word is evaluated at the indicated pair of maps.
Maps belonging to distinct factors have disjoint supports
and commute, so $\rho_\eps$ is a homomorphism.
Let $p_j:\F^k\to\F$ be the coordinate projections, and put
$$
 \lambda_j=6\mu(D_{j,1}\setminus D_{j,2})
              \mu(D_{j,1}\cap D_{j,2})
              \mu(D_{j,2}\setminus D_{j,1})>0.
$$
The factor six counts the assignments of the three labeled
points to the three regions.

The configuration calculation of
\cite[Lemmas~15 and~21]{Brandenbursky2026} gives
\begin{equation}\label{eq:rp2-local-limit}
 \rho_\eps^*\overline\Gamma_{b,3}^{\,q}(v)
    \longrightarrow
       \sum_{j=1}^k\lambda_jp_j^*\jmath_1^*v
       \qquad(\eps\to0)
\end{equation}
in reduced norm. Indeed, outside the transition collars
the endpoints fix every selected point. When the three
points occupy the three displayed regions of one pair,
their trajectory homomorphism is $\jmath_1\circ p_j$ up to
conjugacy. In every other region its image is amenable:
the relevant full twists involve nested or disjoint sets
of points and commute, as in the cited proofs.

For $1\le j\le k$, let $\tau_j:\F^k\to\F^k$ replace
coordinate $j$ by the identity, and put
$$
 \Pi_k=\prod_{j=1}^k(\id-\tau_j^*).
$$
Equation~\eqref{eq:rp2-local-limit} and continuity of cup
products give
\begin{equation}\label{eq:rp2-product-limit}
 \Pi_k\rho_\eps^*\Phi_{1,m,q}(v)
 \longrightarrow
 \Bigl(\prod_{j=1}^{m+1}\lambda_j\Bigr)
 \Sym_{m+1}\bigl(\eta_1^{\times m}\times\jmath_1^*v\bigr).
\end{equation}
The operator $\Pi_k$ kills every term missing a coordinate.
The surviving terms use each coordinate exactly once and
give the displayed symmetrization, with positive signs
because the fixed factors have even degree; see
\cite[Section~4.1]{Brandenbursky2026}.
By Proposition~\ref{prop:external-products}, applied with
$\theta=\eta_1$, the map
$w\mapsto\Sym_{m+1}(\eta_1^{\times m}\times w)$ is injective.
Since $\jmath_1^*$ is injective on $V_q$ and every $\lambda_j$
is positive, \eqref{eq:rp2-product-limit} shows that
$\Phi_{1,m,q}$ is injective.

By Lemma~\ref{lem:rp2-source},
$H^q(\widehat B_3(N_1);\R)=0$ for every $q\ge2$.
Hence every bounded cohomology class of this group
in these degrees is exact.
Proposition~\ref{prop:smooth-transfer} implies
that the restriction of $\overline\Gamma_{b,3}^{\,q}(v)$
to $\Diff_0(N_1,\omega)$ is reduced exact.
Its cup product with the restriction of
$A_1^{\smile m}$ is also reduced exact.
Both $V_2$ and $V_3$ are infinite dimensional, and
the degrees $2m+2$ and $2m+3$, for $m\ge0$, exhaust
all integers at least two.
This proves both assertions for $N_1$.
\end{proof}

\begin{proof}[Proof of Theorem~\ref{thm:bounded-main} for $N_2$]
Use the embedding $\jmath_2:\F\to\widehat B_2(N_2)$ from
Lemma~\ref{lem:klein-source}.
For $q=2,3$, choose an infinite-dimensional linear subspace
$V_q\subset\rEH^q(\widehat B_2(N_2))$ on which $\jmath_2^*$
is injective, and put
\begin{equation}\label{eq:klein-two-point-transfer}
 \Theta_{2,q}=\overline\Gamma_{b,2}^{\,q}|_{V_q}:
 V_q\longrightarrow\rHb^q(\mathcal T_2).
\end{equation}
We first construct homomorphisms $s_\delta:\F\to
\Diff_0(N_2,\omega)$ and a constant $\kappa_2>0$ such that,
for every $\xi\in V_q$, $q=2,3$,
$$
 s_\delta^*\Theta_{2,q}(\xi)
       \longrightarrow\kappa_2\jmath_2^*\xi
$$
in reduced norm. Their defining isotopies will be supported
outside a fixed disc $D_0$. The construction follows the two-point argument of
\cite[Lemma~5.3]{Kimura2025} and \cite[Lemma~26]{Brandenbursky2026}.

Choose simple loops representing $u$ and $v$, as in
Lemma~\ref{lem:klein-source}, which meet only at their
common base point $z_2$ and whose complement in $N_2$
is an open disc containing $z_1$.
Fix a sufficiently small closed disc $U$ about $z_2$.
Choose smooth regular neighborhoods $W_\delta$ of
the union of $U$ and these two loops, with the bands
outside $U$ of width tending to zero.
Each $W_\delta$ is homeomorphic to $N_{2,1}$, contains
$U$ in its interior, and has a complementary disc
containing $z_1$.
Choose a closed disc $D_0$ about $z_1$ in all these
complements.
There are open discs $U_\delta\Subset\operatorname{int}U$
and $V_\delta\Subset N_2\setminus W_\delta$, with
$D_0\subset V_\delta$, such that
$$
 \begin{gathered}
 \mu(U_\delta)\longrightarrow\mu(U)\in(0,1),
 \qquad \mu(V_\delta)\longrightarrow1-\mu(U),\\
 b_\delta:=\mu\bigl(N_2\setminus
                  (U_\delta\cup V_\delta)\bigr)
       \longrightarrow0.
 \end{gathered}
$$
The region outside the two discs consists of the
shrinking bands and boundary collars; see
Figure~\ref{fig:klein-two-point-neighborhood}.

\begin{figure}[tbp]
\centering
\begin{tikzpicture}[x=1.48cm,y=1.48cm,
  font=\small,line cap=round,line join=round,
  >={Stealth[length=1.8mm,width=1.3mm]}]
\def\kleinW{%
 (-0.20,3)--(0.20,3)--(0.20,1.25)
 .. controls (0.20,0.82) and (0.82,0.20) .. (1.25,0.20)
 --(3,0.20)--(3,-0.20)--(1.25,-0.20)
 .. controls (0.82,-0.20) and (0.20,-0.82) .. (0.20,-1.25)
 --(0.20,-3)--(-0.20,-3)--(-0.20,-1.25)
 .. controls (-0.20,-0.82) and (-0.82,-0.20) .. (-1.25,-0.20)
 --(-3,-0.20)--(-3,0.20)--(-1.25,0.20)
 .. controls (-0.82,0.20) and (-0.20,0.82) .. (-0.20,1.25)
 --cycle}
\fill[black!9] \kleinW;
\foreach \angle in {0,90,180,270}{
 \begin{scope}[rotate=\angle]
  \draw[line width=0.65pt] (0.20,3)--(0.20,1.25)
   .. controls (0.20,0.82) and (0.82,0.20) .. (1.25,0.20)
   --(3,0.20);
  \draw[densely dotted,line width=0.85pt] (0.36,3)--(0.36,1.33)
   .. controls (0.36,0.94) and (0.94,0.36) .. (1.33,0.36)
   --(3,0.36);
 \end{scope}
}
\draw[black!65,line width=0.55pt] (0,-3)--(0,3);
\draw[black!65,line width=0.55pt] (-3,0)--(3,0);
\draw[black!65,line width=0.55pt,->] (0,1.83)--(0,2.14);
\draw[black!65,line width=0.55pt,->] (1.83,0)--(2.14,0);
\node[anchor=west,fill=white,inner sep=1.2pt] at (0.50,2.22) {$u$};
\node[anchor=south,fill=white,inner sep=1.2pt] at (2.24,0.49) {$v$};
\draw[dashed,line width=0.65pt] (0,0) circle[radius=0.66];
\draw[densely dotted,line width=0.85pt] (0,0) circle[radius=0.49];
\node[anchor=center,fill=white,inner sep=0.8pt] at (0,-0.28) {$U_\delta$};
\fill (0,0) circle[radius=0.042];
\node[anchor=south west,fill=white,inner sep=0.8pt] at (0.055,0.055) {$z_2$};
\draw[line width=0.55pt] (1.83,1.57) ellipse[x radius=0.35,y radius=0.43];
\fill (1.83,1.57) circle[radius=0.042];
\node[anchor=west,inner sep=1.5pt] at (1.88,1.57) {$z_1$};
\node[anchor=south west] at (2.04,1.95) {$D_0$};
\node at (-1.90,1.70) {$V_\delta$};
\node at (1.43,2.47) {$V_\delta$};
\node at (-1.90,-1.70) {$V_\delta$};
\node at (1.90,-1.70) {$V_\delta$};
\node[anchor=east] (Ulabel) at (-1.12,-0.83) {$U$};
\draw[line width=0.4pt,->] (Ulabel.east) to[out=20,in=215] (-0.48,-0.45);
\node[anchor=east] (Wlabel) at (-0.95,0.99) {$W_\delta$};
\draw[line width=0.4pt,->] (Wlabel.east) to[out=0,in=160] (-0.13,1.37);
\draw[line width=0.8pt] (-3,-3) rectangle (3,3);
\draw[line width=0.8pt,->] (-1.62,3)--(-1.12,3);
\draw[line width=0.8pt,->] (-1.62,-3)--(-1.12,-3);
\draw[line width=0.8pt,->] (-3,1.10)--(-3,1.49);
\draw[line width=0.8pt,->] (-3,1.30)--(-3,1.69);
\draw[line width=0.8pt,->] (3,1.69)--(3,1.30);
\draw[line width=0.8pt,->] (3,1.49)--(3,1.10);
\node[anchor=south] at (0,3.15) {$N_2$};
\draw[line width=0.65pt] (-2.70,-3.43)--(-2.10,-3.43);
\node[anchor=west] at (-2.00,-3.43) {$\partial W_\delta$};
\draw[dashed,line width=0.65pt] (-0.60,-3.43)--(0,-3.43);
\node[anchor=west] at (0.10,-3.43) {$\partial U$};
\draw[densely dotted,line width=0.85pt] (1.15,-3.43)--(1.75,-3.43);
\node[anchor=west] at (1.85,-3.43) {$\partial U_\delta,\partial V_\delta$};
\end{tikzpicture}
\caption{The regions in the Klein-bottle construction. Arrows indicate
side identifications. The shaded neighborhood $W_\delta$ contains the
fixed disc $U$ and the smaller disc $U_\delta$, centered at $z_2$.
The four portions labeled $V_\delta$ form one disc after identification;
this disc contains $D_0$ and $z_1$ and is disjoint from $W_\delta$.
The region outside $U_\delta\cup V_\delta$ consists of the thin bands
and collars and has measure $b_\delta$.
The loops $u$ and $v$ preserve and reverse local orientation,
respectively. The free generators satisfy
$\jmath_2(a)=\overline{v^2}$ and
$\jmath_2(b)=\overline{vuv^{-1}u^{-1}}$.}
\label{fig:klein-two-point-neighborhood}
\end{figure}
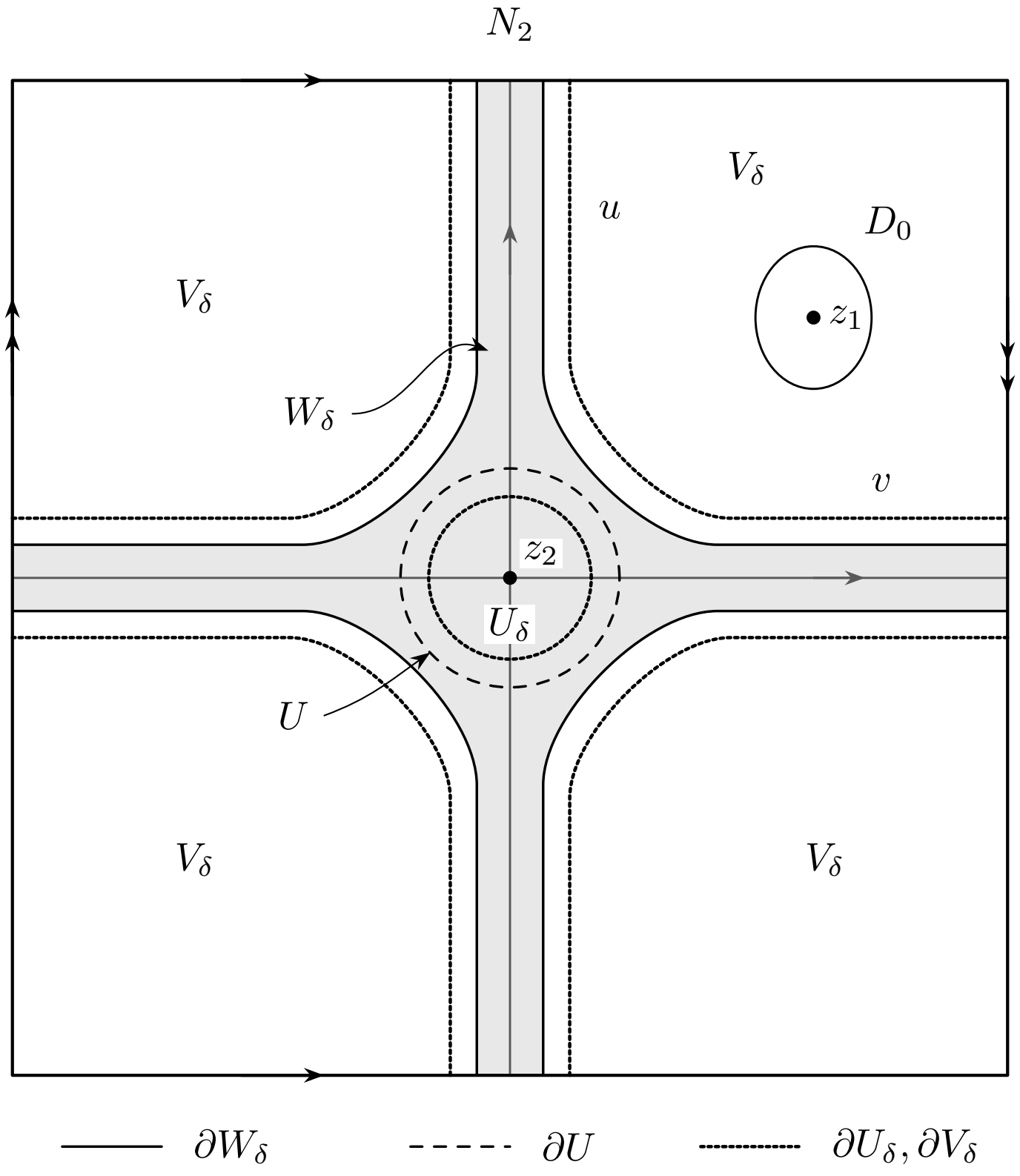

The two loops in \eqref{eq:klein-free-lifts} preserve
local orientation.
For each of the generators $a,b$ of $\F$, push the parametrized disc
$U$ around the inverse of its displayed loop, using an
ordinary smooth isotopy supported in $W_\delta$.
First shrink $U$, move the smaller parametrized disc along the loop,
returning with its original parametrization, and reverse the shrinking.
Use the same initial shrinking for both generators.
The isotopy extension theorem gives an ambient isotopy
whose endpoint is the identity near $U$ and whose entire
path is the identity near $\partial W_\delta$.

We make these maps density preserving by the relative
Moser argument on the orientation cover
$$
 \pi_\delta:\widetilde W_\delta\longrightarrow W_\delta.
$$
Let $\tau$ be its deck involution and let $f$ be one
of the preceding endpoints.
The densities $\omega$ and $f^*\omega$ agree near
$U$ and $\partial W_\delta$, and have equal total mass
on $W_\delta\setminus\operatorname{int}U$.
This complement is connected and nonorientable, so
its orientation cover is connected as well.
Represent the lifted densities by positive area forms
$\Omega_0,\Omega_1$ on $\widetilde W_\delta$.
Both satisfy $\tau^*\Omega_i=-\Omega_i$; their difference has zero
integral and vanishes near $\pi_\delta^{-1}(U)$ and the boundary.
Choose a compactly supported primitive $\theta$ there,
and put
$$
 \theta^-=\tfrac12(\theta-\tau^*\theta),\qquad
 \Omega_r=(1-r)\Omega_0+r\Omega_1,\qquad
 \iota_{X_r}\Omega_r=-\theta^-.
$$
Both $\Omega_r$ and $\theta^-$ change sign under $\tau$,
so $X_r$ is invariant under $\tau$.
Its flow descends to an isotopy $h_r$ on $W_\delta$
which fixes a neighborhood of $U$ and of the boundary,
and satisfies $h_1^*f^*\omega=\omega$.
Extend it by the identity to $N_2$.
This is the relative version of the construction in
\cite[Lemma~6 and proof of Theorem~7]{BruverisMichorParusinskiRainer2018},
with the primitive supported in the indicated complement.

Append $f h_r$ to the original isotopy and denote the
resulting smoothly reparametrized path by $F_t$.
Its endpoints preserve $\omega$.
Apply the same Moser construction to $F_t^*\omega$,
now allowing the correction throughout
$\operatorname{int}W_\delta$.
Choose the primitives using the continuous linear operator of the cited
Lemma~6 before antisymmetrizing them.
The resulting family $H_{t,r}$ is smooth and satisfies
$$
 H_{t,1}^*F_t^*\omega=\omega,\qquad
 H_{t,0}=\id,\qquad H_{0,r}=H_{1,r}=\id.
$$
Thus $F_tH_{t,1}$ is a density-preserving isotopy with
the same endpoints as $F_t$.
The family $F_tH_{t,r}$ gives a homotopy relative to
the endpoints, supported in $W_\delta$.
The first correction fixes $U$ throughout, and the
second preserves the trajectory braids by this homotopy.
Consequently, on configurations in
$X_2(U_\delta\cup V_\delta)$, the trajectory braids
are those of the original disc pushes.

Send $a,b$ to the two corrected endpoints and extend by
evaluating words in these generators. This defines $s_\delta$.
The defining isotopies are supported outside $D_0$,
and every endpoint fixes $U_\delta\cup V_\delta$
pointwise.
For a configuration with one point in each disc, the
homomorphism
$$
 \F\longrightarrow\widehat B_2(N_2),\qquad w\longmapsto\gamma_2(s_\delta(w),x)
$$
is conjugate to $\jmath_2$ after identifying the point
in $V_\delta$ with $z_1$ and the point in $U_\delta$
with $z_2$.
Interchanging the labels conjugates it in the full braid
group and gives the same cohomology contribution.
If both points lie in $V_\delta$, this homomorphism
is trivial.

Suppose that both points lie in $U_\delta$.
Fix such a configuration $x=(x_1,x_2)$ and use it as
the base configuration in this paragraph.
Let $T\in P_2(N_2)$ be the full twist of these two
points inside $U_\delta$.
Every endpoint $s_\delta(w)$ fixes $U_\delta$ pointwise.
Transporting the loop representing $T$ along the
defining isotopy of $s_\delta(w)$ therefore returns
the same loop at the endpoint. The pure trajectory braid
defined by this isotopy therefore commutes with $T$.
Thus these pure braids lie in the centralizer
$$
 C=\{\beta\in P_2(N_2):\beta T=T\beta\}.
$$

We show that $C$ is amenable.
Forgetting the second strand gives the exact sequence
$$
 1\longrightarrow
 F:=\pi_1(N_2\setminus\{x_1\},x_2)
 \longrightarrow P_2(N_2)
 \xrightarrow{\pi}\pi_1(N_2,x_1)
 \longrightarrow1,
$$
where $F$ is free of rank two; see
\cite[Proposition~3.1]{GoncalvesGuaschiLaass2020}.
The full twist $T$ is a nontrivial element of $F$:
it is represented by moving the second point once
around the first.
Hence $C\cap F$ is cyclic, being the centralizer
of a nontrivial element in a free group.
The exact sequence
$$
1\longrightarrow C\cap F\longrightarrow C\longrightarrow\pi(C)\longrightarrow1
$$
proves that $C$ is amenable, since $\pi(C)$ is a subgroup of the
virtually abelian Klein-bottle group. In fact, $\pi(C)$ and $C$
are abelian, but this stronger fact is not needed here. 

Consequently, the image of the trajectory homomorphism for $x$
in $\widehat B_2(N_2)$ lies in the amenable image of $C$.
Changing the base configuration or the closing paths
conjugates this homomorphism and does not affect
the conclusion.
As in the cited two-point calculation, closing paths
inside $U_\delta$ identify these homomorphisms up to
a common conjugation. The bounded change-of-section
homotopy gives the same cohomology classes for the
fixed closing paths.
Thus this region contributes zero in positive-degree
bounded cohomology.

The mixed regions have total measure
$$
\kappa_{2,\delta}=2\mu(U_\delta)\mu(V_\delta)\longrightarrow\kappa_2:=2\mu(U)(1-\mu(U))>0.
$$
The remaining configurations have at least one point
outside $U_\delta\cup V_\delta$, and their measure
is at most $2b_\delta$.
This set is invariant under the endpoints, so its
restricted transfer has norm at most $2b_\delta$.
The preceding calculation gives
$$
\bigl\|s_\delta^*\Theta_{2,q}(\xi)-\kappa_{2,\delta}\jmath_2^*\xi\bigr\|\le 2b_\delta\|\xi\|.
$$
In particular, setting
$e_\delta=2b_\delta+|\kappa_{2,\delta}-\kappa_2|$,
we obtain
\begin{equation}\label{eq:klein-two-point-limit}
\bigl\|s_\delta^*\Theta_{2,q}(\xi)-\kappa_2\jmath_2^*\xi\bigr\|\le e_\delta\|\xi\|,
\qquad e_\delta\longrightarrow0.
\end{equation}
The maps and constants are independent of $q$ and $\xi$.
Since $\jmath_2^*$ is injective on $V_q$ and $\kappa_2>0$,
this proves that $\Theta_{2,q}$ is injective for $q=2,3$.

We next explain exactness of the two-point classes.
The source spaces $V_q$ consist of reduced exact bounded
classes. Proposition~\ref{prop:smooth-transfer} therefore
shows that the classes $\Theta_{2,q}(\xi)$ restrict to reduced exact
classes on $\Diff_0(N_2,\omega)$ for $q=2,3$.
For homeomorphisms, let $\xi\in V_2$ and choose a homogeneous
quasimorphism $\varphi$ on $\widehat B_2(N_2)$ with
$\xi=[\delta\varphi]$.
By Lemma~\ref{lem:klein-source}\textup{(i)}, $\varphi$ vanishes on
the disc braid subgroup.
Theorem~\ref{thm:homeo-main} therefore gives the associated
Gambaudo--Ghys quasimorphism.
Before homogenization, it is
$$
\Psi_\varphi(f)=\int_{X_2(N_2)}\varphi\bigl(\gamma_2(f,x)\bigr)\,d\mu^2(x).
$$
The cocycle relation and invariance of $\mu^2$ show that
its bounded coboundary represents $\Theta_{2,2}(\xi)$.
Thus $\Theta_{2,2}$ has reduced exact image also on
$\Homeo_0(N_2,\mu)$.

To obtain the higher-degree classes, fix $m\ge1$.
Choose $\alpha_2$ and $\eta_2\ne0$ as in
Lemma~\ref{lem:klein-source}\textup{(ii)}, and let
$A_2=\overline\Gamma_{b,3}^{\,2}(\alpha_2)\in\rHb^2(\mathcal T_2)$.
For $q=2,3$, define
\begin{equation}\label{eq:klein-mixed-map}
 \Phi_{2,m,q}:V_q\longrightarrow\rHb^{2m+q}(\mathcal T_2),\qquad
 \Phi_{2,m,q}(\xi)=A_2^{\smile m}\smile\Theta_{2,q}(\xi).
\end{equation}
By Fubini, this is the transfer of
$\alpha_2^{\times m}\times\xi$ over
$X_3(N_2)^m\times X_2(N_2)$, as in
\cite[Section~5.2]{Brandenbursky2026}.

Choose pairwise disjoint discs
$D_1,\ldots,D_m\Subset\operatorname{int}D_0$.
In each $D_j$, choose overlapping round discs
$D_{j,1},D_{j,2}$ and the rotations $f_{j,1,\eps},f_{j,2,\eps}$
as in the proof for $N_1$. Define
$$
 \rho_{j,\eps}:\F\longrightarrow\Diff_0(N_2,\omega),\qquad
 \rho_{j,\eps}(w)=w(f_{j,1,\eps},f_{j,2,\eps}).
$$
The images of the $\rho_{j,\eps}$ have pairwise disjoint supports,
also disjoint from the support of the image of $s_\delta$. Hence
\begin{equation}\label{eq:klein-product-homomorphism}
 \begin{aligned}
 \rho_{\eps,\delta}&:\F^{m+1}
                     \longrightarrow\Diff_0(N_2,\omega),\\
 \rho_{\eps,\delta}(w_1,\ldots,w_m,h)
   &=\rho_{1,\eps}(w_1)\cdots\rho_{m,\eps}(w_m)s_\delta(h)
 \end{aligned}
\end{equation}
is a homomorphism. Put
$\lambda_j=6\mu(D_{j,1}\setminus D_{j,2})\mu(D_{j,1}\cap D_{j,2})\mu(D_{j,2}\setminus D_{j,1})>0$.
The same local calculation as in
\eqref{eq:rp2-local-limit}, now using $\jmath_D$, gives
\begin{equation}\label{eq:klein-disc-limit}
 \rho_{j,\eps}^*A_2\longrightarrow\lambda_j\eta_2
       \qquad(\eps\to0)
\end{equation}
in reduced norm; see
\cite[Lemma~21]{Brandenbursky2026}.

Let $p_j:\F^{m+1}\to\F$, $1\le j\le m+1$, be the
coordinate projections.
For a fixed $\xi\in V_q$, put
$\eta_{j,\eps}=\rho_{j,\eps}^*A_2,\quad E_\delta=s_\delta^*A_2,\quad\xi_\delta=s_\delta^*\Theta_{2,q}(\xi)$.
Thus $\eta_{j,\eps}\in\rHb^2(\F)$,
$E_\delta\in\rHb^2(\F)$, and
$\xi_\delta\in\rHb^q(\F)$.
The degree-two splitting of
\cite[Lemma~22]{Brandenbursky2026} gives
\begin{equation}\label{eq:klein-degree-two-splitting}
 \rho_{\eps,\delta}^*A_2
    =\sum_{j=1}^m p_j^*\eta_{j,\eps}
                   +p_{m+1}^*E_\delta.
\end{equation}
Moreover, the proof of
\cite[Lemma~27]{Brandenbursky2026} gives
\begin{equation}\label{eq:klein-two-point-factorization}
 \rho_{\eps,\delta}^*\Theta_{2,q}(\xi)
       =p_{m+1}^*\xi_\delta.
\end{equation}
Indeed, partition $X_2(N_2)$ according to whether both
points lie in one $D_j$. On the complement, at most one
point lies in each supporting disc, and the local
one-point trajectories are contractible. Choosing the
local parts of the closing paths in these discs removes
their contribution. If both points lie in one $D_j$,
$s_\delta$ is stationary and the disc contribution takes
values in a conjugate of $P_2(D_j)\cong\Z$.
Its bounded-cohomology contribution vanishes by amenability.
The change of closing paths and the integration of the
bounded primitives are exactly those in the cited proof.

For $1\le j\le m$, let $\tau_j:\F^{m+1}\to\F^{m+1}$
replace coordinate $j$ by the identity, and put
$$
 \Pi_m=\prod_{j=1}^m(\id-\tau_j^*).
$$
Expand \eqref{eq:klein-mixed-map} using
\eqref{eq:klein-degree-two-splitting} and
\eqref{eq:klein-two-point-factorization}.
There are $m$ degree-two factors to distribute.
A term using $E_\delta$, or using one disc coordinate twice,
misses another disc coordinate and is killed by $\Pi_m$.
The remaining $m!$ terms use each disc coordinate exactly
once, and their signs agree because the factors have even
degree. Thus
\begin{equation}\label{eq:klein-mixed-limit}
\Pi_m\rho_{\eps,\delta}^*\Phi_{2,m,q}(\xi)=m!\,\eta_{1,\eps}\times\cdots\times\eta_{m,\eps}
\times\xi_\delta\longrightarrow m!\kappa_2\Bigl(\prod_{j=1}^m\lambda_j\Bigr)
\eta_2^{\times m}\times\jmath_2^*\xi.
\end{equation}
The limit is taken in reduced norm as $(\eps,\delta)\to(0,0)$, by
\eqref{eq:klein-two-point-limit} and
\eqref{eq:klein-disc-limit}.
This is the mixed-transfer calculation of
\cite[Section~5.1, equations~(68)--(70), and Section~5.2,
equations~(79)--(80)]{Brandenbursky2026}.

By Proposition~\ref{prop:external-products}, applied with
$\theta=\eta_2$, external product with $\eta_2^{\times m}$
is injective on $\rHb^q(\F)$.
The map $\jmath_2^*$ is injective on $V_q$, and the coefficient
in \eqref{eq:klein-mixed-limit} is positive.
It follows that $\Phi_{2,m,q}(\xi)=0$ implies $\xi=0$.

The two-point factor is reduced exact for diffeomorphisms
when $q=2,3$, and for homeomorphisms when $q=2$, as proved
above. Its cup product with $A_2^{\smile m}$ is 
reduced exact in these cases. Explicitly, if a bounded
cocycle $d_\xi=\delta B_\xi$ represents an exact lift of
this factor and $c_A$ represents $A_2$, then
$$
c_A^{\smile m}\smile d_\xi=\delta\bigl(c_A^{\smile m}\smile B_\xi\bigr).
$$
Thus exactness holds before reduction.
The infinite-dimensional spaces $V_2$ and $V_3$, together with
the injectivity of $\Theta_{2,q}$ and $\Phi_{2,m,q}$, give
the required reduced bounded cohomology in every degree,
reduced exact bounded cohomology in every degree for
diffeomorphisms, and reduced exact bounded cohomology
in every positive even degree for homeomorphisms.
\end{proof}

We conclude with the following question.
\begin{question*}
For every odd integer $d>3$, is the reduced exact bounded cohomology
$\rEH^d\bigl(\Homeo_0(N_2,\mu)\bigr)$ infinite-dimensional?
\end{question*}

\bibliographystyle{plain}
\bibliography{references}

\begin{thebibliography}{10}

\bibitem{BBF2016}
Mladen Bestvina, Ken Bromberg, and Koji Fujiwara.
\newblock Stable commutator length on mapping class groups.
\newblock {\em Ann. Inst. Fourier (Grenoble)}, 66(3):871--898, 2016.

\bibitem{BestvinaFujiwara2002}
Mladen Bestvina and Koji Fujiwara.
\newblock Bounded cohomology of subgroups of mapping class groups.
\newblock {\em Geom. Topol.}, 6:69--89, 2002.

\bibitem{BestvinaFujiwara2007}
Mladen Bestvina and Koji Fujiwara.
\newblock Quasi-homomorphisms on mapping class groups.
\newblock {\em Glasnik Mat. Ser. III}, 42(62)(1):213--236, 2007.

\bibitem{Boke2025}
Lukas B{\"o}ke.
\newblock The {Klein} bottle has stably unbounded homeomorphism group.
\newblock {\em Glasgow Math. J.}, 67(3):494--499, 2025.

\bibitem{Brandenbursky2015}
Michael Brandenbursky.
\newblock Bi-invariant metrics and quasi-morphisms on groups of {Hamiltonian} diffeomorphisms of surfaces.
\newblock {\em Internat. J. Math.}, 26(9):1550066, 2015.
\newblock 29 pp.

\bibitem{Brandenbursky2026}
Michael Brandenbursky.
\newblock Bounded cohomology of transformation groups in all degrees, 2026.
\newblock Preprint, \href{https://arxiv.org/abs/2610.11327}{arXiv:2610.11327}.

\bibitem{BrandenburskyNonorientable2026}
Michael Brandenbursky.
\newblock Exact bounded cohomology of measure-preserving homeomorphism groups of nonorientable hyperbolic manifolds, 2026.
\newblock Preprint, \href{https://arxiv.org/abs/2610.32475}{arXiv:2610.32475}.

\bibitem{BrandenburskyKedra2013}
Michael Brandenbursky and Jarek K{\k{e}}dra.
\newblock On the autonomous metric on the group of area-preserving diffeomorphisms of the 2-disc.
\newblock {\em Algebr. Geom. Topol.}, 13(2):795--816, 2013.

\bibitem{BKMS2026}
Michael Brandenbursky, Jarek K{\k e}dra, Micha{\l} Marcinkowski, and Egor Shelukhin.
\newblock On quasimorphisms and distortion in homeomorphism groups.
\newblock {\em Commun. Contemp. Math.}, page 2650050, 2026.

\bibitem{Entropy}
Michael Brandenbursky and Micha{\l} Marcinkowski.
\newblock Entropy and quasimorphisms.
\newblock {\em J. Mod. Dyn.}, 15:143--163, 2019.

\bibitem{BrandenburskyMarcinkowski2022}
Michael Brandenbursky and Micha{\l} Marcinkowski.
\newblock Bounded cohomology of transformation groups.
\newblock {\em Math. Ann.}, 382(3--4):1181--1197, 2022.

\bibitem{BrandenburskyMenashe2026}
Michael Brandenbursky and Lior Menashe.
\newblock Bounded cohomology of measure-preserving homeomorphism groups of non-orientable surfaces, 2026.
\newblock Preprint, \href{https://arxiv.org/abs/2506.02728}{arXiv:2506.02728v3}.

\bibitem{Brooks1981}
Robert Brooks.
\newblock Some remarks on bounded cohomology.
\newblock In {\em Riemann Surfaces and Related Topics: Proceedings of the 1978 Stony Brook Conference}, volume~97 of {\em Annals of Mathematics Studies}, pages 53--63. Princeton University Press, Princeton, NJ, 1981.

\bibitem{BruverisMichorParusinskiRainer2018}
Martins Bruveris, Peter~W. Michor, Adam Parusi{\'n}ski, and Armin Rainer.
\newblock Moser's theorem on manifolds with corners.
\newblock {\em Proc. Amer. Math. Soc.}, 146(11):4889--4897, 2018.

\bibitem{Calegari2009}
Danny Calegari.
\newblock {\em {scl}}, volume~20 of {\em MSJ Memoirs}.
\newblock Mathematical Society of Japan, Tokyo, 2009.

\bibitem{CalegariFujiwara2010}
Danny Calegari and Koji Fujiwara.
\newblock Stable commutator length in word-hyperbolic groups.
\newblock {\em Groups Geom. Dyn.}, 4(1):59--90, 2010.

\bibitem{EntovPolterovich2003}
Michael Entov and Leonid Polterovich.
\newblock {Calabi} quasimorphism and quantum homology.
\newblock {\em Int. Math. Res. Not.}, 2003(30):1635--1676, 2003.

\bibitem{FarbLubotzkyMinsky2001}
Benson Farb, Alexander Lubotzky, and Yair Minsky.
\newblock Rank-1 phenomena for mapping class groups.
\newblock {\em Duke Math. J.}, 106(3):581--597, 2001.

\bibitem{Farre2018}
James Farr{\'e}.
\newblock Bounded cohomology of finitely generated {Kleinian} groups.
\newblock {\em Geom. Funct. Anal.}, 28(6):1597--1640, 2018.

\bibitem{Fathi1980}
Albert Fathi.
\newblock Structure of the group of homeomorphisms preserving a good measure on a compact manifold.
\newblock {\em Ann. Sci. {\'E}cole Norm. Sup. (4)}, 13(1):45--93, 1980.

\bibitem{FranksHandel2006}
John Franks and Michael Handel.
\newblock Distortion elements in group actions on surfaces.
\newblock {\em Duke Math. J.}, 131(3):441--468, 2006.

\bibitem{Frigerio}
Roberto Frigerio.
\newblock {\em Bounded Cohomology of Discrete Groups}, volume 227 of {\em Mathematical Surveys and Monographs}.
\newblock American Mathematical Society, Providence, RI, 2017.

\bibitem{FrigerioPozzettiSisto2015}
Roberto Frigerio, Maria~Beatrice Pozzetti, and Alessandro Sisto.
\newblock Extending higher-dimensional quasi-cocycles.
\newblock {\em J. Topol.}, 8(4):1123--1155, 2015.

\bibitem{GambaudoGhys2004}
Jean-Marc Gambaudo and {\'E}tienne Ghys.
\newblock Commutators and diffeomorphisms of surfaces.
\newblock {\em Ergodic Theory Dynam. Systems}, 24(5):1591--1617, 2004.

\bibitem{GoncalvesGuaschi2004}
Daciberg~Lima Gon{\c c}alves and John Guaschi.
\newblock The braid groups of the projective plane.
\newblock {\em Algebr. Geom. Topol.}, 4:757--780, 2004.

\bibitem{GoncalvesGuaschiLaass2020}
Daciberg~Lima Gon{\c c}alves, John Guaschi, and Vinicius~Casteluber Laass.
\newblock The {Borsuk--Ulam} property for homotopy classes of maps from the torus to the {Klein} bottle.
\newblock {\em Topol. Methods Nonlinear Anal.}, 56(2):529--558, 2020.

\bibitem{GoncalvesGuaschiMaldonado2018}
Daciberg~Lima Gon{\c c}alves, John Guaschi, and Miguel Maldonado.
\newblock Embeddings and the (virtual) cohomological dimension of the braid and mapping class groups of surfaces.
\newblock {\em Confluentes Math.}, 10(1):41--61, 2018.

\bibitem{GuaschiPereiro2018}
John Guaschi and Carolina de~Miranda~e Pereiro.
\newblock Lower central and derived series of semi-direct products, and applications to surface braid groups.
\newblock {\em J. Pure Appl. Algebra}, 224(7):106309, 2020.

\bibitem{GuaschiJuanPineda2026}
John Guaschi and Daniel Juan-Pineda.
\newblock Braid groups of the projective plane, mapping class groups of non-orientable surfaces, and algebraic {$K$}-theory of their group rings.
\newblock {\em Bol. Soc. Mat. Mex.}, 32:Paper No. 56, 2026.

\bibitem{HullKapovich2018}
Michael Hull and Ilya Kapovich.
\newblock Counting conjugacy classes in {$\operatorname{Out}(F_N)$}.
\newblock {\em Bull. Aust. Math. Soc.}, 97(3):412--421, 2018.

\bibitem{Ishida2014}
Tomohiko Ishida.
\newblock Quasi-morphisms on the group of area-preserving diffeomorphisms of the 2-disk via braid groups.
\newblock {\em Proc. Amer. Math. Soc. Ser. B}, 1:43--51, 2014.

\bibitem{KatayamaKuno2024}
Takuya Katayama and Erika Kuno.
\newblock The mapping class group of a nonorientable surface is quasi-isometrically embedded in the mapping class group of the orientation double cover.
\newblock {\em Groups Geom. Dyn.}, 18(2):407--418, 2024.

\bibitem{KimMaruyamaFlux2025}
KyeongRo Kim and Shuhei Maruyama.
\newblock The {Euler} class and flux homomorphisms under non-orientability, 2025.
\newblock Preprint, \href{https://arxiv.org/abs/2508.12874}{arXiv:2508.12874}.

\bibitem{KimMaruyamaIMRN2026}
KyeongRo Kim and Shuhei Maruyama.
\newblock Quasimorphisms on the group of density preserving diffeomorphisms of the {M\"obius} band.
\newblock {\em Int. Math. Res. Not. IMRN}, 2026(1):rnaf365, 2026.

\bibitem{Kimura2025}
Mitsuaki Kimura.
\newblock {Gambaudo--Ghys} construction on bounded cohomology.
\newblock {\em J. Math. Soc. Japan}, 77(1):135--152, 2025.

\bibitem{Kumar2021}
Rishi Kumar.
\newblock Construction of quasi-morphisms on groups of measure-preserving homeomorphisms of non-orientable surfaces.
\newblock Master's thesis, Ben Gurion University of the Negev, 2021.

\bibitem{MatsumotoMorita1985}
Shigenori Matsumoto and Shigeyuki Morita.
\newblock Bounded cohomology of certain groups of homeomorphisms.
\newblock {\em Proc. Amer. Math. Soc.}, 94(3):539--544, 1985.

\bibitem{McDuffSalamon2017}
Dusa McDuff and Dietmar Salamon.
\newblock {\em Introduction to Symplectic Topology}, volume~27 of {\em Oxford Graduate Texts in Mathematics}.
\newblock Oxford University Press, Oxford, third edition, 2017.

\bibitem{Nitsche}
Martin Nitsche.
\newblock Higher-degree bounded cohomology of transformation groups, 2021.
\newblock arXiv:2105.08698.

\bibitem{OxtobyUlam1941}
John~C. Oxtoby and Stanislaw~M. Ulam.
\newblock Measure-preserving homeomorphisms and metrical transitivity.
\newblock {\em Ann. of Math. (2)}, 42(4):874--920, 1941.

\bibitem{ParisRolfsen1999}
Luis Paris and Dale Rolfsen.
\newblock Geometric subgroups of surface braid groups.
\newblock {\em Ann. Inst. Fourier (Grenoble)}, 49(2):417--472, 1999.

\bibitem{Polterovich2006}
Leonid Polterovich.
\newblock {Floer} homology, dynamics and groups.
\newblock In Paul Biran, Octav Cornea, and Fran{\c{c}}ois Lalonde, editors, {\em Morse theoretic methods in nonlinear analysis and in symplectic topology}, volume 217 of {\em NATO Sci. Ser. II Math. Phys. Chem.}, pages 417--438. Springer, Dordrecht, 2006.

\bibitem{Shtern2001}
Alexander~I. Shtern.
\newblock Remarks on pseudocharacters and the real continuous bounded cohomology of connected locally compact groups.
\newblock {\em Ann. Global Anal. Geom.}, 20(3):199--221, 2001.

\bibitem{Soma1997}
Teruhiko Soma.
\newblock Bounded cohomology and topologically tame {Kleinian} groups.
\newblock {\em Duke Math. J.}, 88(2):357--370, 1997.

\bibitem{Stukow2017}
Micha{\l} Stukow.
\newblock Subgroups generated by two {Dehn} twists on a nonorientable surface.
\newblock {\em Topology Proc.}, 50:151--201, 2017.

\end{thebibliography}
\end{document}